\documentclass[10pt,a4paper]{article}
\usepackage{geometry}
\usepackage[utf8]{inputenc}
\usepackage{amsthm,amssymb,amsmath,natbib,graphicx,caption}
\usepackage{xcolor,enumitem}
\usepackage[ruled,vlined]{algorithm2e}
\SetAlFnt{\small}\SetAlCapFnt{\small}
\setcitestyle{numbers,comma,square}
\usepackage{tikz}
\usetikzlibrary{matrix,calc}
\usepackage[pdfpagelabels,hypertexnames=false,bookmarksopen=true,
bookmarksopenlevel=1,colorlinks=true,pdfstartview=FitV,
linkcolor=black,citecolor=black,urlcolor=black]{hyperref}
\usepackage{cleveref}
\hypersetup{pdftitle={A Deterministic Algorithm for Eigenvalue Separation in Tridiagonal Matrices},pdfauthor={Wenqi Zhu and Ping Nie}}

\newtheorem{proposition}{Proposition}[section]
\newtheorem{theorem}{Theorem}[section]

\newtheorem{lemma}{Lemma}[section]
\newtheorem{example}{Example}[section]
\newtheorem{remark}{Remark}[section]

\newtheorem{definition}{Definition}[section]
\newcommand{\norm}[1]{\left\lVert #1\right\rVert}
\DeclareMathOperator{\gap}{gap}
\DeclareMathOperator{\diag}{diag}
\DeclareMathOperator{\tr}{tr}

\DeclareMathOperator{\spn}{span}
\newcommand{\R}{\mathbb R}
\newcommand{\C}{\mathbb C}

\title{A Deterministic Algorithm for Eigenvalue Separation\\in Tridiagonal Matrices}
\author{Wenqi Zhu\thanks{Corresponding author. Mathematical Institute, University of Oxford, Andrew Wiles Building, Woodstock Road, Oxford OX2 6GG, United Kingdom. \texttt{wenqi.zhu@maths.ox.ac.uk}} \and
Ping Nie\thanks{David R. Cheriton School of Computer Science, University of Waterloo, 200 University Avenue West, Waterloo, ON N2L 3G1, Canada. \texttt{ping.nie@uwaterloo.ca}}}
\date{November 2026}
\usepackage{xurl}
\usepackage{float}
\usepackage{needspace}
\begin{document}
\maketitle

\begin{abstract}
We develop a deterministic algorithm for separating the eigenvalues
of Hermitian tridiagonal matrices by small diagonal perturbations. 
For $\|T\|_2\leq1$, $0<\delta<1/2$, and
$X=\diag(1,\ldots,n)$, our algorithm selects $t\in[0,\delta/n]$ such
that $D=tX$ satisfies
\[
\|D\|_2\leq\delta,
\qquad
\gap(T+D)>c \delta^2n^{-12}
\]
where $c>0$ is a constant.
The construction uses a classical convex weighted eigenvalue sum and its
finite-scale slope variation. A tridiagonal localization argument turns
near-degenerate eigenvalues into separated responses to the fixed diagonal
direction, yielding a quantitative gap certificate. An adaptive interval
selector then identifies a certifying window in at most
$\lceil\log_2(16n^3)\rceil$ refinements, using $O(\log n)$
evaluations of the weighted eigenvalue sum to inverse-polynomial
accuracy. A structure-exploiting implementation uses rigorous
Sturm eigenvalue enclosures and adaptive precision to certify the
stopping condition on six families of test matrices. The directly
certified gaps are substantially larger than the worst-case analytical
bounds, illustrating both the effectiveness of the construction and the
conservatism of the universal certificate. Our results give a
deterministic answer to the eigenvalue-separation question underlying
the deterministic Minami perturbation problem for tridiagonal
matrices.
\end{abstract}

\section{Introduction}

\label{sec:problem}

Spectral separation plays an important role in numerical linear
algebra. When eigenvalues are tightly clustered, computing their
individual eigenvectors and associated invariant subspaces can
become sensitive to perturbations. Spectral gaps also influence
the accuracy, convergence, and deflation criteria of many
eigenvalue algorithms \cite{Parlett1996,Nakatsukasa2012,
HorningTownsend2020}. A natural approach is to perturb the
matrix slightly so that its eigenvalues become better separated.
The challenge is to construct such a perturbation efficiently,
without first computing an expensive eigendecomposition.

Diagonal perturbations are particularly attractive for
tridiagonal matrices because they preserve the matrix structure
and the cost of matrix-vector products. In random spectral
theory, Minami-type estimates provide probabilistic bounds on
the occurrence of closely clustered eigenvalues under suitable
independent diagonal perturbations
\cite{Minami1996,CombesGerminetKlein2009}. Whether comparable
separation can be achieved deterministically is a more
challenging question. A closely related question,
\emph{deterministic pseudospectral shattering}, a deterministic
counterpart of the random-perturbation shattering of Banks et
al.~\cite{BanksEtAl2023}, was raised by
Sobczyk~\cite{Sobczyk2025} (Question~4) in his deterministic
complexity analysis of Hermitian eigenproblems; the present paper
answers its eigenvalue-gap version for Hermitian tridiagonal
matrices.  The question was
subsequently posed as Problem~3.2,
\emph{Deterministic Minami Perturbation}, in the report of the
Simons Institute workshop \emph{Linear Systems and Eigenvalue
Problems}, held in Berkeley in October 2025
\cite{AmselEtAl2026}. The KE-02 entry in the Open Problems
in NLA catalogue further asks for a nearly linear deterministic
construction for arbitrary normalized Hermitian tridiagonal
matrices \cite{CatalogueKE02}.

Let $n\geq2$, $0<\delta<1/2$, and let $T\in\C^{n\times n}$
be Hermitian tridiagonal with $\norm{T}_2\leq1$. Thus
\begin{equation}\label{eq:T}
T=\begin{pmatrix}
d_1&\beta_1&&&\\
\overline{\beta_1}&d_2&\beta_2&&\\
&\overline{\beta_2}&d_3&\ddots&\\
&&\ddots&\ddots&\beta_{n-1}\\
&&&\overline{\beta_{n-1}}&d_n
\end{pmatrix},\qquad d_i\in\R.
\end{equation}
We write its eigenvalues in nondecreasing order and
define its minimum eigenvalue gap by
\begin{equation}\label{eq:gap}
\gap(T):=\min_{1\leq j<n}
\bigl(\lambda_{j+1}(T)-\lambda_j(T)\bigr).
\end{equation}
In particular, $\gap(T)=0$ if $T$ has a repeated eigenvalue.
The objective is to construct a real diagonal matrix $D$
such that
\begin{equation}\label{eq:target}
\norm{D}_2\leq\delta,\qquad
\gap(T+D)\geq c_0(\delta/n)^p,
\end{equation}
for absolute constants $c_0,p>0$, at a total cost of
$O\bigl(n[1+\log(n/\delta)]^q\bigr)$ exact arithmetic
operations for a constant $q$.
This cost concerns the construction of $D$, including
any internal spectral evaluations; it does not require
outputting an eigenvector matrix. Arithmetic complexity
and finite-precision complexity are separate questions.

Several related approaches provide insight into this problem.
Classical work by Ye~\cite{Ye1995} and
Parlett~\cite{Parlett1996} explains how closely spaced
eigenvalues of tridiagonal matrices are associated with
the localization of their invariant subspaces.
Nakatsukasa~\cite{Nakatsukasa2012} uses bounds on
eigenvector components to study the sensitivity of
eigenvalues to structured perturbations.
More recently, deterministic constructions have been
obtained for special classes, including weakly coupled
tridiagonal matrices and matrices with constant diagonal
and equal off-diagonal magnitudes
\cite{CatalogueKE02}. These results provide important
structural insights, but do not establish the requested
construction for the unrestricted
tridiagonal case.

In this paper, we give a deterministic construction for
the full Hermitian tridiagonal class. Our approach takes
a different perspective: rather than estimating individual
eigenvalue movements directly, we use the convexity of a
weighted spectral function to detect and avoid regions of
strong spectral sensitivity.
In particular, we restrict the perturbation to a
single predetermined diagonal direction:
\begin{equation}\label{eq:tilt}
\boxed{
X=\diag(1,\ldots,n),\qquad
H(t)=T+tX,\qquad 0\leq t\leq\delta/n.
}
\end{equation}
Since $\norm{tX}_2=nt\leq\delta$, the search reduces to selecting
one scalar $t$ at which spectral separation can be certified without
first computing the minimum gap accurately. 

\paragraph{Main results and contributions.}
This paper delivers two main results.

\emph{First, we develop a convex spectral-gap certificate.}
We consider a weighted sum of the ordered eigenvalues
of $H(t)$, viewed as a function $F(t)$ of the perturbation
parameter. This spectral function is convex, strictly
increasing, and has bounded slopes with respect to the perturbation size $t$. Up to an affine
trace term, it measures the total pairwise spectral
distance after perturbation. By the Hellmann--Feynman formula, its first
derivative is a weighted sum of the responses of
individual eigenvalues to the diagonal perturbation.
Moreover, classical second-order perturbation formulas
show that the curvature contains pairwise contributions
of the form
\[
\underbrace{\text{curvature contribution}}_{
\text{for a pair of eigenvalues}}
\ \propto\
\frac{|\text{perturbation coupling}|^2}
{\text{eigenvalue gap}}.
\]
Thus closely spaced eigenvalues can produce a strong
change in spectral slope, provided that their responses
to the perturbation are sufficiently different.

The difficulty is that a small gap alone does not
guarantee large curvature: the corresponding coupling
may be weak or vanish. To overcome this, we exploit
the tridiagonal structure. A discrete Wronskian
identity and a small-leakage cut convert a nearly
degenerate eigenvalue pair into two spatially
separated approximate eigenvectors. Their responses
to the fixed diagonal direction $X$ are quantitatively
different. Using a nearby auxiliary matrix and the
variational representation of the spectral function,
we turn this response separation into a finite-scale
gap certificate.

This proof combines ideas from classical eigenvector
localization and structured perturbation theory
\cite{Ye1995,Parlett1996,Nakatsukasa2012}
with convex spectral analysis
\cite{OvertonWomersley1993,Lewis1996Convex,
LewisSendov2001}. The key is to use these classical
tools for a different purpose: instead of bounding
the sensitivity of a given eigenvalue or optimizing
a spectral objective, we certify a lower bound on
the minimum spectral gap.

\emph{Second, we construct an adaptive deterministic
interval-selection algorithm.}
Starting from the admissible parameter interval
$[0,\delta/n]$, the algorithm repeatedly selects
one of two shorter windows, each retaining slightly
less than half of the current interval. Convexity
allows us to compare their endpoint slope variations
and choose a window with a smaller spectral response.
At every iteration, we directly test whether the
centered slope variation is sufficiently small to
certify spectral separation. Once this condition
is satisfied, the algorithm terminates. 

Each refinement requires only a constant number
of evaluations of the spectral function. In
particular, the selection rule does not require
prior knowledge of the eigenvectors or a full
eigendecomposition at each step. We prove that
the algorithm terminates after at most
\[
K=\left\lceil\log_2(16n^3)\right\rceil
=O(\log n)
\]
refinements. The returned perturbation satisfies
\begin{equation}\label{eq:intro-main-result}
\boxed{
\norm{tX}_2\leq\delta,\qquad
\gap(T+tX)\geq c\,\frac{\delta^2}{n^{12}},
}
\end{equation}
where $c>0$ is a constant. 
The construction uses at most $3+8K$ evaluations of the spectral
function, each to an absolute accuracy no finer than
$\delta/(2^{16}Kn^4)$. Its cost is therefore that of $O(\log n)$
computations of all eigenvalues of a Hermitian tridiagonal matrix to
inverse-polynomial accuracy (\Cref{subsec:oracle}). 
The two contributions are complementary:
the first provides a local certificate relating
spectral response to eigenvalue separation, while
the second uses convexity to find a certifying
perturbation parameter efficiently.

\paragraph{Numerical validation.}
We implement the adaptive algorithm using double-precision window
selection, rigorous Sturm-based eigenvalue enclosures, and adaptive
precision. On six families of tridiagonal matrices the certified
stopping condition is typically satisfied well before the theoretical
iteration bound. Matrices with exactly repeated or extremely close
eigenvalues are separated to certified gaps above $10^{-5}$, and the
directly certified gaps are many orders of magnitude larger than the
analytical lower bounds, decaying only mildly with $n$ over the tested
range. These results demonstrate the practical effectiveness of the
construction on the tested instances, without implying an improved
worst-case exponent.

\paragraph{Related work: localization and structured perturbation.}
Classical work on close eigenvalues of tridiagonal
matrices provides an important foundation for our
geometric argument. Ye~\cite{Ye1995} relates
closely spaced eigenvalues to complementary
principal submatrices, even when the off-diagonal
entry connecting the two blocks is not small.
Parlett~\cite{Parlett1996} studies the localization
of invariant subspaces associated with tightly
clustered eigenvalues, while V{\"o}mel and
Parlett~\cite{VoemelParlett2011} investigate
how such localization can be detected.

These results show that small spectral gaps
are closely connected to the spatial structure
of the corresponding eigenspaces. Our localization
argument follows this viewpoint, but requires a
specific quantitative conclusion: two approximate
eigenvectors supported on opposite sides of
a coordinate cut, with small residuals and
separated responses to the predetermined
direction $X$. 
Nakatsukasa~\cite{Nakatsukasa2012} develops
a complementary perturbation approach for
Hermitian block tridiagonal matrices. By
bounding eigenvector components through the
block recurrence, he obtains upper bounds
on the movement of individual eigenvalues
under structured perturbations, with applications
to Wilkinson matrices and aggressive early
deflation. In contrast, we use the localization
of a near-degenerate pair to establish a lower
bound on its collective spectral response.
This response is then converted into a
minimum-gap certificate through convexity.
Thus our work connects the geometric ideas
of tridiagonal localization to a different
spectral-certification problem.

\paragraph{Related work: convex eigenvalue optimization.}
Weighted eigenvalue sums and diagonal matrix
corrections have a long history in spectral
optimization, including graph partitioning
and nondifferentiable optimization
\cite{DonathHoffman1973,CullumDonathWolfe1975,
LewisOverton1996}.
Ky Fan's variational principle and subsequent
work of Overton and Womersley, Lewis, and
Lewis--Sendov establish the convexity,
subdifferential, and second-order calculus
of spectral functions
\cite{Fan1949,OvertonWomersley1992,
OvertonWomersley1993,Lewis1996Convex,
Lewis1996Derivatives,LewisSendov2001}. 
These classical results are key to our spectral function. In eigenvalue
optimization, such functions have typically
been used to characterize optimality, derive
subgradients, or construct optimization
algorithms. Here we use the same variational
structure in a different direction. 
In particular, our proof uses the existence
of distinct supporting spectral slopes at
a nearby matrix with a repeated eigenvalue.
The difference between these slopes is
controlled by the responses of two localized
eigenvectors to $X$. We then transfer this
information back to the original matrix
using stability of the convex spectral
function. This is the central connection
between structured eigenvalue perturbation
and nonsmooth convex analysis in our proof. 

Rather than minimizing the spectral function,
we search for a window on which its finite-scale
slope variation is small. Such a window certifies
that the spectrum is sufficiently separated.
The interval-selection mechanism is also related
to adaptive sampling of convex functions.
Simchowitz et al.~\cite{SimchowitzEtAl2018},
for example, use midpoint secant errors to
guide the approximation of a convex function.
Their midpoint error is directly related to
our centered slope variation. Related spectral and Gram-certificate characterizations of
regularized polynomial subproblems, and variational methods for tensor
eigenvalues, appear in our earlier work
\cite{CartisZhu2026MathComp,CartisZhu2026SoSExact,ChenZhuCartis2026Dinkelbach,ZhuChenZhou2026HTM,Zhu2026Thesis}.

\paragraph{Notation.}
We use $\norm{\cdot}_2$ for the Euclidean vector norm
and the matrix operator norm, $\norm{\cdot}_{\mathrm F}$
for the Frobenius norm, and $\norm{\cdot}_*$ for the
nuclear norm, the sum of the singular values.
The vector $e_i$ is the $i$th standard unit vector.
At every permitted parameter,
$\norm{H(t)}_2\leq1+\delta<2$, and its off-diagonal
entries have absolute value at most one.

\begin{remark}[Reduction to real symmetric matrices]
\label{rem:real}
A diagonal unitary matrix $S$ can be chosen so that
$S^*TS$ is real symmetric. Choose its phases
successively along each chain of nonzero
off-diagonal entries; zero entries separate
the chains. Since $S^*XS=X$, the entire family
$T+tX$ is preserved up to unitary similarity.
Thus the eigenvalues, gaps, and norms are
unchanged, and the real symmetric arguments
below also apply to Hermitian tridiagonal matrices.
\end{remark}

The remainder of the paper is organized as follows.
Section~\ref{sec:convex-spectral} develops the
convex spectral function and its finite-scale
slope variation. Section~\ref{sec:local-gap}
uses tridiagonal localization to establish
the local gap certificate.
Section~\ref{sec:selection} develops an adaptive
window selector and analyzes spectral-function
evaluations.
Section~\ref{sec:numerics} presents rigorously
certified numerical experiments.
Section~\ref{sec:conclusion} summarizes
the results and discusses possible extensions.
\section{The convex spectral function and its slope variation}
\label{sec:convex-spectral}

Throughout this section, the results apply to any Hermitian matrix, not only to tridiagonal matrices. We use the notation $T$ for simplicity. 

\paragraph{The spectral function.}
For the fixed family $H(t)=T+tX$, define the spectral function
\begin{equation}\label{eq:F}
\boxed{F(t):=\sum_{j=1}^n j\lambda_j(H(t)).}
\end{equation}
The map $A\mapsto\sum_{j=1}^n j\lambda_j(A)$ is a convex
spectral function, meaning that it is convex in a Hermitian
matrix $A$ and depends only on its eigenvalues.
Its basic properties follow from Ky Fan's maximum principle
\cite{Fan1949}, which provides a variational characterization
of sums of the largest eigenvalues. Overton and Womersley
\cite{OvertonWomersley1992} further studied these eigenvalue
sums and their convex-analytic properties, while Lewis
\cite{Lewis1996Convex} developed a general convex-analytic
framework for spectral functions on Hermitian matrices.

These variational characterizations have played an
important role in eigenvalue optimization. In particular,
Overton and Womersley~\cite{OvertonWomersley1993}
developed optimality conditions and duality theory
for minimizing sums of the largest eigenvalues.
Closely related weighted trace objectives also appear
in Brockett's work on matrix gradient flows,
eigenvalue ordering, and diagonalization
\cite{Brockett1991}. The weights $1,\ldots,n$
in our spectral function correspond to such a
weighted trace objective through \eqref{eq:variational}.

While these spectral functions and their convexity
properties are classical, our interest is different:
we study their variation along a fixed diagonal
perturbation and use this variation to certify
quantitative eigenvalue separation.

\paragraph{Why this spectral function?}
A useful interpretation follows by collecting coefficients:
\begin{equation}\label{eq:pairwise-spread}
F(t)=\frac{n+1}{2}\tr H(t)
+\frac12\sum_{i<j}\bigl(\lambda_j(t)-\lambda_i(t)\bigr).
\end{equation}
Thus $F$ measures the total pairwise spectral distance,
up to an affine trace term. \Cref{ex:tails} illustrates its
decomposition into sums of the largest eigenvalues and shows
why we instead use its change under perturbation to control
the smallest gap.

\begin{example}[Tail sums and spectral spread]\label{ex:tails}
When $n=3$, the tail-sum decomposition in \eqref{eq:variational} is
\[
F(t)=\lambda_1(t)+2\lambda_2(t)+3\lambda_3(t)
=\underbrace{\lambda_1(t)+\lambda_2(t)+\lambda_3(t)}_{\tr H(t)}
+\underbrace{\lambda_2(t)+\lambda_3(t)}_{\text{two largest}}
+\underbrace{\lambda_3(t)}_{\text{largest}}.
\]
The significance of $F$ lies in its variation under the perturbation $tX$.
As $t$ varies, the eigenvalues move at different rates,
and their pairwise gaps change accordingly. When two
eigenvalues approach one another, their relative motion
may change rapidly, producing curvature in $F$.
This motivates studying the second difference of $F$,
which removes the affine trace contribution and measures
the variation of the spectral spread over a finite
perturbation window. We will use this variation to
derive a certificate for the minimum eigenvalue gap.
\end{example}

The following lemma collects the technical properties of the spectral function needed for our analysis. Note that the variational principle and convexity in \Cref{lem:convex} are classical
\cite{Fan1949,OvertonWomersley1992,Lewis1996Convex}. We give the proof for completeness. 

\begin{lemma}[Properties of the spectral function]\label{lem:convex}
For every Hermitian $T$, the function $F$ is
\textbf{strictly increasing} and \textbf{convex}.
For $u<v$, its secant slopes satisfy
\begin{equation}\label{eq:slopes}
\frac{n(n+1)}2
\leq\frac{F(v)-F(u)}{v-u}
\leq\frac{n^2(n+1)}2.
\end{equation}
Moreover, it admits the variational representation
\begin{equation}\label{eq:variational}
F(t)
=\max_{U^*U=I}\tr\bigl(UXU^*(T+tX)\bigr)
=\tr(T+tX)+\sum_{j=2}^n\sum_{\ell=j}^n\lambda_\ell(t).
\end{equation}
\end{lemma}

\begin{proof}
Write $\lambda_\ell(t)=\lambda_\ell(T+tX)$. Counting the
coefficient of each eigenvalue gives
$
F(t)=\tr(T+tX)+\sum_{j=2}^n\sum_{\ell=j}^n\lambda_\ell(t).
$
Ky Fan's maximum principle \cite{Fan1949} states that
\[
\sum_{\ell=n-k+1}^n\lambda_\ell(t)
=\max_{\substack{V\in\C^{n\times k}\\V^*V=I_k}}
\tr\bigl(V^*(T+tX)V\bigr).
\]
For fixed $V$, the expression on the right is affine in $t$.
Its pointwise maximum is therefore convex, and \eqref{eq:variational}
proves the convexity of $F$. To obtain the first equality
in \eqref{eq:variational}, expand the weighted trace in the
same way into a trace term and tail sums of Rayleigh quotients.
Each tail sum is bounded by the corresponding eigenvalue sum,
and all bounds are attained simultaneously in an ordered
orthonormal eigenbasis. When $T$ is real symmetric, the
maximizing matrix may be chosen real orthogonal. 
Finally, since $I\preceq X\preceq nI$, each ordered eigenvalue
increases by an amount between $v-u$ and $n(v-u)$ for $u<v$.
Thus $F$ is strictly increasing with bounded secant slopes.
\end{proof}

\paragraph{Spectral sensitivity under perturbation.}

The connection between the spectral function and
eigenvalue separation can also be understood through
classical perturbation theory. Suppose that $H(t)$
has simple eigenvalues, with normalized eigenvectors
$u_j(t)$. The Hellmann--Feynman formula gives
\begin{equation}\label{eq:spectral-first-derivative}
\lambda_j'(t)=u_j(t)^*Xu_j(t),\qquad
F'(t)=\sum_{j=1}^n j\,u_j(t)^*Xu_j(t).
\end{equation}
Thus the derivative of $F$ is a weighted sum of the
responses of individual eigenvalues to the diagonal
perturbation. Differentiating again and pairing the
eigenvalue terms yields
\begin{equation}\label{eq:spectral-second-derivative}
F''(t)=2\sum_{j<k}
\frac{(k-j)|u_j(t)^*Xu_k(t)|^2}
{\lambda_k(t)-\lambda_j(t)}.
\end{equation}
Schematically, each nonzero pairwise contribution obeys
\[
\underbrace{\text{curvature contribution}}_{\text{for a pair of eigenvalues}}
\ \propto\
\frac{|\text{perturbation coupling}|^2}{\text{eigenvalue gap}},
\]
with the precise weight $2(k-j)$ in \eqref{eq:spectral-second-derivative}.
The numerator measures the coupling between the two
eigenvectors induced by $X$, while the denominator
is their eigenvalue separation. Thus small gaps can
produce large curvature when the corresponding
perturbation coupling is non-negligible. Conversely,
a small gap need not produce large curvature if this
coupling is weak. These observations are related to
the differential theory of spectral functions
\cite{Lewis1996Derivatives,LewisSendov2001,
DrusvyatskiyPaquette2018}.

This motivates a finite-scale measure of spectral
sensitivity that remains meaningful when eigenvalues
are repeated and classical derivatives are unavailable.

\paragraph{Computing and interpreting slope variation.}
\begin{definition}[Centered slope variation]\label{def:defect}
For $h>0$, define the centered slope variation of $F$ at $t$ by
\begin{equation}\label{eq:defect}
\boxed{
\mathcal S_h(t)
:=
\frac{F(t+h)-F(t)}h
-
\frac{F(t)-F(t-h)}h
\geq0.
}
\end{equation}
\end{definition}
The quantity $\mathcal S_h(t)$ requires three
spectral-function evaluations and measures the difference
between the right and left secant slopes over the window
$[t-h,t+h]$.  

Equivalently, the slope variation can be written
directly in terms of the ordered eigenvalues:
\begin{equation}\label{eq:slope-eigenvalues}
\mathcal S_h(t)=\sum_{j=1}^n j
\left(
\frac{\lambda_j(t+h)-\lambda_j(t)}{h}
-\frac{\lambda_j(t)-\lambda_j(t-h)}{h}
\right).
\end{equation}
Hence $\mathcal S_h(t)$ measures the change in the
weighted average eigenvalue slopes between the
left and right halves of the window. When the
eigenvalues respond approximately linearly to the
perturbation, this variation is small. Rapid changes
in their relative motion, including changes in the
ordering of continuously tracked eigenvalue branches,
can instead produce a large slope variation.

If $F$ is twice continuously differentiable on the window,
the slope variation can also be expressed as
\begin{equation}\label{eq:slope-curvature}
\mathcal S_h(t)
=
\int_{-h}^{h}
\left(1-\frac{|s|}{h}\right)F''(t+s)\,ds.
\end{equation}
Thus $\mathcal S_h(t)$ is a weighted
integral of the curvature over the prescribed window.

The perturbation scale offers another interpretation.
If the family is extended to large positive $t$, the
term $tX$ dominates the fixed matrix $T$, and the
eigenvalues satisfy $\lambda_j(t)=jt+O(1)$.
In fact, $F$ becomes asymptotically affine, so
$\mathcal S_h(t)\to0$ for each fixed $h>0$ as
$t\to+\infty$. This is a useful limiting picture,
but it lies outside the small-perturbation regime
required by our problem.

Here $0\leq t\leq\delta/n$, so domination by
$tX$ cannot be assumed. Instead, the convexity
and bounded total slope variation of $F$ allow us
to search for a window with sufficiently stable
spectral response. The local certificate in
\Cref{prop:bridge} shows that, for our normalized
tridiagonal family,
\[
\mathcal S_h(t)<\frac12
\quad\Longrightarrow\quad
\gap(H(t))>\frac{h^2}{2^{20}n^4}.
\]
Thus small slope variation certifies quantitative
separation at the window center. The procedure
in \Cref{sec:selection} then finds such a window.

\begin{figure}[h]
\centering
\includegraphics[width=0.98\linewidth]{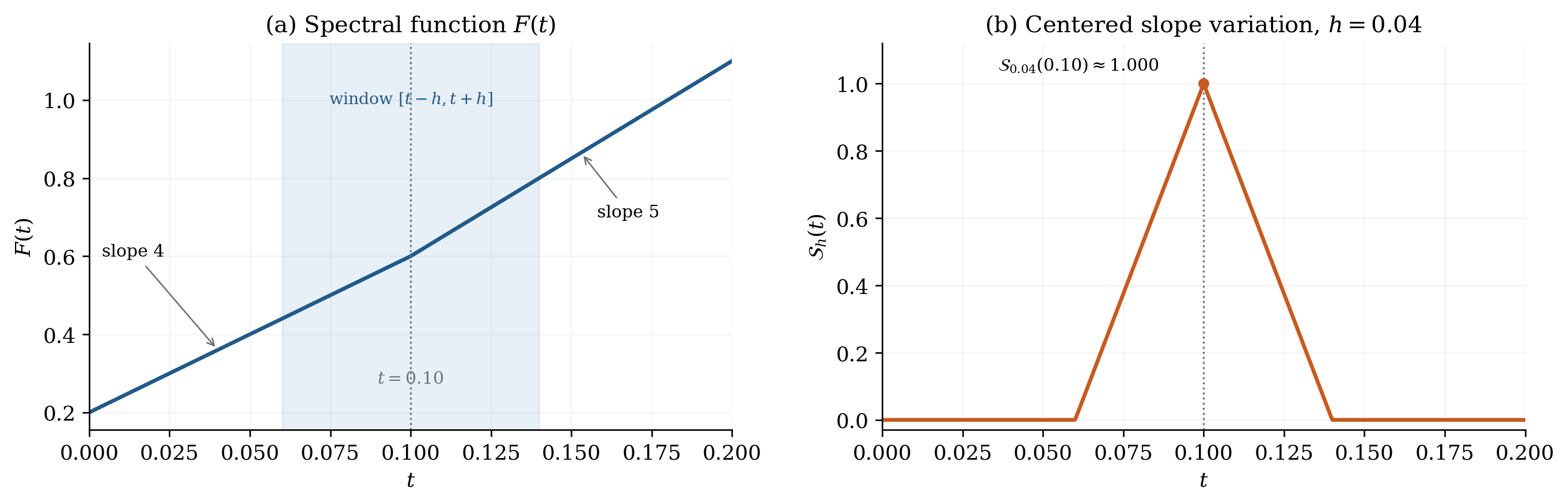}
\caption{\Cref{ex:collision} with $h=0.04$.
Left: the spectral function $F(t)$, with the window
$[t-h,t+h]$ highlighted around $t=0.1$.
The slope increases from $4$ to $5$ at the center.
Right: the corresponding centered slope variation
$\mathcal S_h(t)$, which forms a triangular peak
of height $1$ at $t=0.1$.}
\label{fig:convexity-defect}
\end{figure}

\begin{example}[A collision produces a corner]
\label{ex:collision}
Let $T=\diag(0.1,0)$ and $X=\diag(1,2)$.
Then $H(t)=\diag(0.1+t,2t)$, whose eigenvalues coincide
at $t=0.1$. Ordering them gives
\begin{equation}\label{eq:collision}
F(t)=
\begin{cases}
0.2+4t,&t\leq0.1,\\
0.1+5t,&t\geq0.1,
\end{cases}
\qquad
F(t)=0.15+4.5t+\frac12|t-0.1|.
\end{equation}
The slope increases from $4$ to $5$. Dividing the centered
second difference by $h$ gives the centered slope variation
\begin{equation}\label{eq:collision-defect-explicit}
\mathcal S_h(t)
=
\frac1{2h}\bigl(
|t+h-0.1|+|t-h-0.1|-2|t-0.1|
\bigr)
=
\max\left\{1-\frac{|t-0.1|}{h},0\right\}.
\end{equation}
Thus $\mathcal S_h$ is a triangular peak centered at $t=0.1$; see \Cref{fig:convexity-defect}. 
Since affine terms cancel in the slope difference, the trace
contribution in \eqref{eq:pairwise-spread} does not affect
$\mathcal S_h(t)$. 
\end{example}

\begin{figure}[h]
\centering
\includegraphics[width=0.98\linewidth]{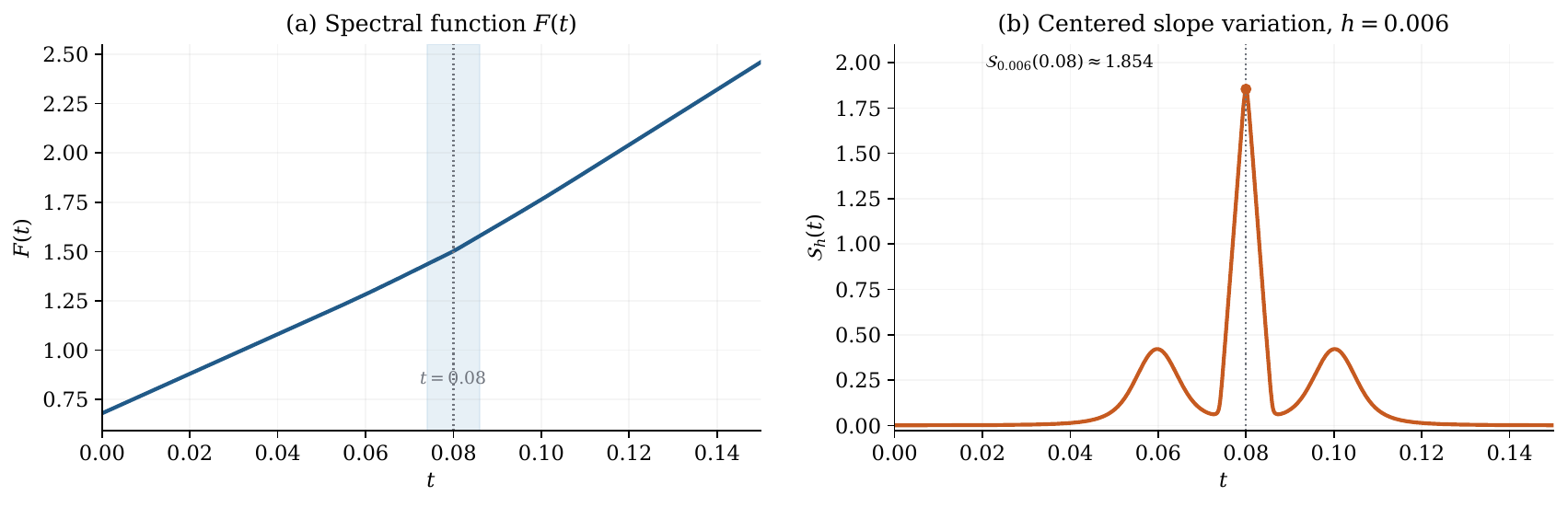}
\caption{\Cref{ex:tridiagonal-defect} with $h=0.006$.
Left: the spectral function $F(t)$, with the
perturbation window highlighted around $t=0.08$.
Right: the corresponding centered slope variation
$\mathcal S_h(t)$, showing three smooth peaks.
At $t=0.08$, the slope variation is approximately
$1.8541$.}
\label{fig:tridiagonal-defect}
\end{figure}

\begin{example}[A tridiagonal matrix with avoided crossings]
\label{ex:tridiagonal-defect}
Let
\begin{equation}\label{eq:tridiagonal-example}
T=
\begin{pmatrix}
0.16&0.003&0\\
0.003&0.10&0.003\\
0&0.003&0
\end{pmatrix},
\qquad
X=\diag(1,2,3).
\end{equation}
Both off-diagonal edges are nonzero, so $T+tX$ has
simple eigenvalues for every $t$. Thus $F$ and $\mathcal S_h(t)$ are smooth.
\Cref{fig:tridiagonal-defect} shows the local bend and
the slope variation over a larger range of centers.
There are three visible peaks near $0.06,0.08,0.10$,
the crossing locations of the diagonal branches when
the edges are set to zero. With the nonzero edges in
\eqref{eq:tridiagonal-example}, these are smooth peaks,
not exact corners. These computations illustrate the
slope variation; they are not certified numerical bounds.
\end{example}

\section{Small slope variation certifies a spectral gap}
\label{sec:local-gap}

The main result of this section shows that a sufficiently small centered
slope variation certifies a quantitative lower bound on the minimum
eigenvalue gap. We first state the certificate and then develop its two
ingredients: a localization argument that uses the tridiagonal structure,
and a stability argument that transfers the slope variation from a nearby
auxiliary matrix to the original perturbation family.

Throughout, fix $t\in[0,\delta/n]$ and write $H=H(t)=T+tX$ when the
parameter is fixed. We work in real symmetric coordinates by
\Cref{rem:real}; in particular, $\norm{H}_2\leq2$ and
$|H_{i,i+1}|\leq1$.

\begin{proposition}[Local gap certificate]\label{prop:bridge}
For $t\in[0,\delta/n]$ and $0<h\leq1$,
\begin{equation}\label{eq:bridge}
\boxed{\mathcal S_h(t)<\frac12\quad\Longrightarrow\quad
\gap(H(t))>\frac{h^2}{2^{20}n^4}.}
\end{equation}
\end{proposition}

We prove the contrapositive. A small eigenvalue gap forces two approximate
eigenvectors to be separated by a cut in the coordinate basis of $X$.
We then make them exact eigenvectors at a common eigenvalue of a nearby
matrix. Their distinct responses to $X$ force a slope variation of at
least one, and the perturbation error is too small to remove this signal.

\subsection{Localization: a near-degenerate pair splits across a cut}
\label{subsec:localization}

The key structural step is to construct two approximate eigenvectors whose
projections lie in disjoint coordinate subspaces. This is where the
tridiagonal assumption enters; note that the tridiagonal assumption is only used in this subsection. 
\Cref{ex:localization} illustrates how two nearby eigenvalues of a
tridiagonal matrix give rise to approximate eigenvectors
supported in disjoint coordinate subspaces, with
distinct responses to the diagonal perturbation.

\begin{example}[An illustrative example of a small-leakage cut]
\label{ex:localization}

Consider an $8\times8$ real symmetric tridiagonal matrix
\[
H=\diag(d_1,\ldots,d_8)
+\sum_{i=1}^7\beta_i
(e_ie_{i+1}^{\mathsf T}+e_{i+1}e_i^{\mathsf T}),
\]
where, with $\gamma=10^{-5}$,
\[
\begin{aligned}
(d_1,\ldots,d_8)
={}&
\left(\frac14-\frac\gamma2,\,
\frac12-\frac\gamma2,\,
\frac14-\frac\gamma2,\,
\frac12,\,\frac12,\,
\frac\gamma2,\,\frac\gamma2,\,\frac\gamma2\right),\\
(\beta_1,\ldots,\beta_7)
={}&
\left(-\frac14,-\frac14,0.0005,0.1,
0.0005,\frac14,\frac14\right).
\end{aligned}
\]
All off-diagonal entries are nonzero, and $\norm{H}_2<1$. 
Two eigenvalues are
$\lambda\approx-5.173792\times10^{-6}$ and
$\mu\approx4.739763\times10^{-6}$, satisfying
$|\lambda-\mu|<\gamma$.
Their corresponding orthonormal eigenvectors are
approximately
\[
\begin{aligned}
u_1&\approx(0.577344,0.577345,0.577346,-0.000602,
0.000123,-0.003033,0,0.003033)^{\mathsf T},\\
u_2&\approx(0.002477,0.002477,0.002477,0.000145,
-0.000736,0.707100,0,-0.707100)^{\mathsf T}.
\end{aligned}
\]

We choose the cut $k=4$, separating the coordinate
subspaces $\spn\{e_1,\ldots,e_4\}$ and
$\spn\{e_5,\ldots,e_8\}$. 
We call $k=4$ a small-leakage cut. 
We obtain two orthonormal approximate eigenvectors
\[
\begin{aligned}
v_1&\approx(0.577350,0.577350,0.577351,-0.000601,
0,0,0,0)^{\mathsf T},\\
v_2&\approx(0,0,0,0,-0.000737,0.707107,
7.37\times10^{-7},-0.707107)^{\mathsf T}.
\end{aligned}
\]
Their coordinate supports are disjoint, and direct
computation gives
\[
\norm{\left(H-\frac{\lambda+\mu}{2}I\right)v_1}_2
\approx6.15\times10^{-5},\qquad
\norm{\left(H-\frac{\lambda+\mu}{2}I\right)v_2}_2
\approx7.53\times10^{-5}.
\]
Moreover, for $X=\diag(1,\ldots,8)$,
\[
v_2^{\mathsf T}Xv_2-v_1^{\mathsf T}Xv_1
\approx5.0000\geq1.
\]
\textbf{How do we obtain $v_1,v_2$?} We first construct two linear
combinations of the original eigenvectors $u_1,u_2$.
Choose $p=2$ and $q=6$, corresponding to two rows
with substantial norms and distinct directions,
on opposite sides of the cut. Write $r_i=((u_1)_i,(u_2)_i)\in\R^{1\times2}$ and define
\[
w_1=\frac{u_1[p]u_1+u_2[p]u_2}
{\sqrt{u_1[p]^2+u_2[p]^2}} = \frac{\gamma_p^T }{\|\gamma_p\|}(u_1, u_2),
\qquad
w_2=\frac{u_1[q]u_1+u_2[q]u_2}
{\sqrt{u_1[q]^2+u_2[q]^2}}= \frac{\gamma_q^T }{\|\gamma_q\|}(u_1, u_2), 
\]
where $u_j[i]$ denotes the $i$th entry of $u_j$.
Let $P_L$ and $P_R$ be the orthogonal projections
onto $\spn\{e_1,\ldots,e_4\}$ and
$\spn\{e_5,\ldots,e_8\}$, respectively. We set
\[
v_1=\frac{P_Lw_1}{\norm{P_Lw_1}_2},
\qquad
v_2=\frac{P_Rw_2}{\norm{P_Rw_2}_2}.
\]
\end{example}

\begin{lemma}[Spatially separated approximate eigenvectors]
\label{lem:localization}
Suppose $\lambda,\mu$ are two eigenvalues of $H(t)$, counted with
multiplicity, satisfying $0\leq|\lambda-\mu|\leq\gamma\leq
1/(256n^2)$ with $\gamma>0$. Then there exist a cut $1\leq k<n$
and orthonormal unit vectors
$v_1\in\spn\{e_1,\ldots,e_k\}$ and
$v_2\in\spn\{e_{k+1},\ldots,e_n\}$ such that
\begin{equation}\label{eq:quasi}
\left\lVert\left(H(t)-\frac{\lambda+\mu}{2}I\right)v_i
\right\rVert_2\leq16n\sqrt\gamma\quad(i=1,2),
\qquad v_2^{\mathsf T}Xv_2-v_1^{\mathsf T}Xv_1\geq1.
\end{equation}
\end{lemma}

\begin{proof}
Fix $t$ and write $H=H(t)$. Choose orthonormal eigenvectors
$u_1,u_2$ associated with $\lambda,\mu$, respectively, and
write $r_i=((u_1)_i,(u_2)_i)\in\R^{1\times2}$.
By orthogonality, $\sum_{i=1}^n r_i^{\mathsf T}r_i=I_2$.

\emph{Step 1. Two substantial rows with separated directions.}
Let $J=\{i:\norm{r_i}_2\geq1/(2\sqrt n)\}$.
The discarded rows have total squared norm at most $1/4$.
Consequently,
$\sum_{i\in J}r_i^{\mathsf T}r_i\succeq(3/4)I_2$
and $\sum_{i\in J}\norm{r_i}_2^2\leq2$.

Choose any $p\in J$ and a unit vector $e\in\R^2$
perpendicular to $r_p$. Then
$\sum_{i\in J}|r_ie|^2\geq3/4$, so there exists
$q\in J$ satisfying
$|r_qe|^2/\norm{r_q}_2^2\geq3/8$.
After exchanging $p,q$ if necessary, assume $p<q$.

For every nonzero row, let
$\Pi_i=r_i^{\mathsf T}r_i/\norm{r_i}_2^2$,
the rank-one orthogonal projector onto its direction.
Since
$\norm{\Pi_i-\Pi_j}_2
=|\det(r_i,r_j)|/(\norm{r_i}_2\norm{r_j}_2)$,
we obtain
\begin{equation}\label{eq:endpoints}
\norm{r_p}_2,\norm{r_q}_2\geq\frac1{2\sqrt n},
\qquad
\norm{\Pi_p-\Pi_q}_2\geq\sqrt{\frac38}>\frac12.
\end{equation}
Thus two rows have substantial norms but distinctly
different directions.

\emph{Step 2. Wronskian identity.}
Write $d_i=H_{ii}$ and $\beta_i=H_{i,i+1}$.
The eigenvector equations at coordinate $i$ are
\[
\begin{aligned}
\beta_{i-1}(u_1)_{i-1}+d_i(u_1)_i+\beta_i(u_1)_{i+1}
&=\lambda(u_1)_i,\\
\beta_{i-1}(u_2)_{i-1}+d_i(u_2)_i+\beta_i(u_2)_{i+1}
&=\mu(u_2)_i.
\end{aligned}
\]
Multiplying the first equation by $(u_2)_i$,
the second by $(u_1)_i$, subtracting, and summing
from $1$ to $i$ gives the discrete Wronskian identity
\begin{equation}\label{eq:wronskian}
\beta_i
\det\begin{pmatrix}
(u_1)_i&(u_1)_{i+1}\\
(u_2)_i&(u_2)_{i+1}
\end{pmatrix}
=
\beta_i\det\begin{pmatrix}r_i\\r_{i+1}\end{pmatrix}
=
(\mu-\lambda)\sum_{j=1}^i(u_1)_j(u_2)_j.
\end{equation}

Let $P_{\leq i}=\sum_{j=1}^i e_je_j^{\mathsf T}$.
Since $u_1^{\mathsf T}u_2=0$, Cauchy--Schwarz gives
\[
\left|\sum_{j=1}^i(u_1)_j(u_2)_j\right|
=
\left|u_1^{\mathsf T}
\left(P_{\leq i}-\frac12I\right)u_2\right|
\leq\frac12.
\]
Thus
$|\beta_i|\,|\det(r_i,r_{i+1})|\leq\gamma/2$.

\emph{Step 3. Finding a small-leakage cut.}
Set $\tau=4\sqrt{n\gamma}$. We claim that there exists
$p\leq k<q$ such that
\begin{equation}\label{eq:leak}
|\beta_k|\max\{\norm{r_k}_2,\norm{r_{k+1}}_2\}
\leq\tau.
\end{equation}
Suppose otherwise. Call a row \emph{large} if
$\norm{r_i}_2>\tau$, and \emph{small} otherwise.
Since $\tau\leq1/(4\sqrt n)$, the endpoint rows
$r_p,r_q$ are large.

We distinguish three cases.

\emph{(a) Two adjacent large rows.}
If $r_i,r_{i+1}$ are both large, the failure of
\eqref{eq:leak} implies
$|\beta_i|\norm{r_i}_2\norm{r_{i+1}}_2>\tau^2$.
The Wronskian bound therefore yields
\[
\norm{\Pi_i-\Pi_{i+1}}_2
=
\frac{|\det(r_i,r_{i+1})|}
{\norm{r_i}_2\norm{r_{i+1}}_2}
\leq\frac{\gamma}{2\tau^2}.
\]
Thus two adjacent large rows must have nearly
parallel directions.

\emph{(b) Two adjacent small rows.}
If both $r_i,r_{i+1}$ are small, then
$|\beta_i|\leq1$ gives
$|\beta_i|\max\{\norm{r_i}_2,\norm{r_{i+1}}_2\}
\leq\tau$, contradicting our assumption.
Hence small rows cannot occur consecutively.

\emph{(c) A large--small--large configuration.}
Suppose $r_{i-1},r_{i+1}$ are large but $r_i$
is small. The direction of $r_i$ itself need
not be well controlled, so we compare the
directions of $r_{i-1}$ and $r_{i+1}$ directly.

The three-term recurrence at row $i$ gives
\[
\beta_{i-1}r_{i-1}
+\left(d_i-\frac{\lambda+\mu}{2}\right)r_i
+\beta_ir_{i+1}
=
r_i\diag\left(\frac{\lambda-\mu}{2},
\frac{\mu-\lambda}{2}\right).
\]
Taking the determinant with $r_{i-1}$,
using the Wronskian bound and
$|d_i-(\lambda+\mu)/2|\leq4$, yields
\[
\norm{\Pi_{i-1}-\Pi_{i+1}}_2
\leq
\frac{2\gamma}
{|\beta_{i-1}\beta_i|
\norm{r_{i-1}}_2\norm{r_{i+1}}_2}
+
\frac{\gamma\norm{r_i}_2}
{2|\beta_i|\norm{r_{i+1}}_2}
\leq\frac{5\gamma}{\tau^2}.
\]
Indeed, the failure of \eqref{eq:leak} implies
$|\beta_{i-1}|\norm{r_{i-1}}_2>\tau$
and $|\beta_i|\norm{r_{i+1}}_2>\tau$,
while $\norm{r_i}_2\leq\tau\leq1$.

We now follow the large rows from $p$ to $q$,
skipping any isolated small rows. By cases (a)--(c),
the total change in direction satisfies
\[
\norm{\Pi_p-\Pi_q}_2
\leq\frac{5n\gamma}{\tau^2}
=\frac5{16}<\frac12.
\]
This contradicts \eqref{eq:endpoints}, proving
the existence of the cut \eqref{eq:leak}.

\emph{Step 4. Projection, truncation, and normalization.}
Using the row $r_p$, define the unit vector
\[
w_1=
\frac{(u_1)_p u_1+(u_2)_p u_2}{\norm{r_p}_2}.
\]
By construction, $(w_1)_p=\norm{r_p}_2$,
$|(w_1)_i|\leq\norm{r_i}_2$, and
$\left\lVert(H-\frac{\lambda+\mu}{2}I)w_1
\right\rVert_2\leq\gamma/2$.

Let $P_L=\sum_{i=1}^k e_ie_i^{\mathsf T}$
be the orthogonal projection onto the first
$k$ coordinate directions. Since $p\leq k$,
$\norm{P_Lw_1}_2\geq\norm{r_p}_2\geq1/(2\sqrt n)$.

By tridiagonality, the commutator
$[H,P_L]=HP_L-P_LH$ has only two boundary
contributions:
\[
[H,P_L]w_1
=\beta_k(w_1)_ke_{k+1}
-\beta_k(w_1)_{k+1}e_k.
\]
The small-leakage bound \eqref{eq:leak}
therefore gives
$\norm{[H,P_L]w_1}_2\leq\sqrt2\tau$.

Define $v_1=P_Lw_1/\norm{P_Lw_1}_2$.
Combining the residual and boundary estimates,
\[
\left\lVert
\left(H-\frac{\lambda+\mu}{2}I\right)v_1
\right\rVert_2
\leq2\sqrt n\left(\frac\gamma2+\sqrt2\tau\right)
\leq4\sqrt n\tau=16n\sqrt\gamma.
\]

Similarly, using the row $r_q$, define
$w_2=((u_1)_q u_1+(u_2)_q u_2)/\norm{r_q}_2$.
Let $P_R=\sum_{i=k+1}^n e_ie_i^{\mathsf T}$
and set $v_2=P_Rw_2/\norm{P_Rw_2}_2$.
The same argument gives
$\left\lVert(H-\frac{\lambda+\mu}{2}I)v_2
\right\rVert_2\leq16n\sqrt\gamma$.

The vectors $v_1,v_2$ have disjoint coordinate
supports and are therefore orthonormal.
Finally, since $X=\diag(1,\ldots,n)$,
$v_1^{\mathsf T}Xv_1\leq k$ and
$v_2^{\mathsf T}Xv_2\geq k+1$, proving
$v_2^{\mathsf T}Xv_2-v_1^{\mathsf T}Xv_1\geq1$.
\end{proof}

\begin{example}[Identifying the small-leakage cut in \Cref{ex:localization}]
\label{rem:localization-example}
Returning to \Cref{ex:localization}, write
$[u_1\ \ u_2]^{\mathsf T}=[\gamma_1\ \cdots\ \gamma_8]$,
where $\gamma_i=((u_1)_i,(u_2)_i)^{\mathsf T}=r_i^{\mathsf T}$.
The determinant in \eqref{eq:wronskian} is the $2\times2$ minor
$\det[\gamma_i\ \gamma_{i+1}]$. The dashed vertical line below
marks the selected cut $k=4$, between columns $4$ and $5$.
\begin{equation}\label{eq:localization-row-minors}
\resizebox{0.91\linewidth}{!}{$\displaystyle
[u_1\ \ u_2]^{\mathsf T}=[\gamma_1\ \cdots\ \gamma_8]\approx
\begin{tikzpicture}[baseline=(M.center)]
\matrix (M) [matrix of math nodes,
 left delimiter={[},right delimiter={]},
 nodes={inner sep=1.6pt,anchor=center},
 column sep=1pt,row sep=3pt,font=\scriptsize] {
 0.577344 \pgfmatrixnextcell 0.577345 \pgfmatrixnextcell 0.577346 \pgfmatrixnextcell -0.000602 \pgfmatrixnextcell 0.000123 \pgfmatrixnextcell -0.003033 \pgfmatrixnextcell -1.23\cdot10^{-7} \pgfmatrixnextcell 0.003033 \\
 0.002477 \pgfmatrixnextcell 0.002477 \pgfmatrixnextcell 0.002477 \pgfmatrixnextcell 0.000145 \pgfmatrixnextcell -0.000736 \pgfmatrixnextcell 0.707100 \pgfmatrixnextcell 7.36\cdot10^{-7} \pgfmatrixnextcell -0.707100 \\
};
\draw[blue!70!black,line width=.75pt,rounded corners=1pt]
 ([xshift=-1.5pt,yshift=2.3pt]M-1-1.north west)
 rectangle ([xshift=1.5pt,yshift=-2.3pt]M-2-2.south east);
\draw[orange!80!black,line width=.75pt,rounded corners=1pt]
 ([xshift=-1.5pt,yshift=2.3pt]M-1-4.north west)
 rectangle ([xshift=1.5pt,yshift=-2.3pt]M-2-5.south east);
\draw[teal!70!black,dashed,line width=.65pt,rounded corners=1pt]
 ([xshift=-1.5pt,yshift=2pt]M-1-6.north west)
 rectangle ([xshift=1.5pt,yshift=-2pt]M-2-7.south east);
\draw[teal!70!black,dashed,line width=.65pt,rounded corners=1pt]
 ([xshift=-1pt,yshift=3.7pt]M-1-7.north west)
 rectangle ([xshift=1pt,yshift=-3.7pt]M-2-8.south east);
\coordinate (cutTop) at ($(M-1-4.north east)!0.5!(M-1-5.north west)$);
\coordinate (cutBottom) at ($(M-2-4.south east)!0.5!(M-2-5.south west)$);
\coordinate (cutBottomAligned) at (cutTop |- M-2-4.south);
\draw[red!80!black,densely dashed,line width=1.25pt]
 ([yshift=7pt]cutTop) -- ([yshift=-7pt]cutBottomAligned);
\node[red!80!black,above=5pt,font=\scriptsize] at (cutTop) {$k=4$};
\end{tikzpicture}.
$}
\end{equation}

For $\gamma=10^{-5}$ and $n=8$, the threshold is
$\tau=4\sqrt{n\gamma}\approx0.0358$. At any cut $k$,
we must check the \emph{leakage}
$L_k:=|\beta_k|\max\{\|\gamma_k\|_2,\|\gamma_{k+1}\|_2\}$,
not just the determinant of the two columns. The values of $L_k$ for all seven cuts are:
\[
\begin{array}{c|ccccccc}
 k&1&2&3&4&5&6&7\\ \hline
 \text{row types}&\mathrm{L-L}&\mathrm{L-L}&\mathrm{L-S}&\mathrm{S-S}&\mathrm{S-L}&\mathrm{L-S}&\mathrm{S-L}\\
 L_k&.1443&.1443&.000289&.0000746&.000354&.1768&.1768\\
 L_k\leq\tau?&\text{No}&\text{No}&\text{Yes}&\text{Yes}&\text{Yes}&\text{No}&\text{No}
\end{array}
\]
Here $\mathrm{L}$ and $\mathrm{S}$ indicate whether the row norm is
larger or smaller than $\tau$, respectively.

\emph{(a) Two adjacent large rows (blue box).}
At $k=1$, both columns have norms approximately $0.577$,
and $|\det[\gamma_1\ \gamma_2]|\approx5.67\times10^{-8}$.
Nevertheless, $L_1\approx0.1443>\tau$, so this is \emph{not}
a small-leakage cut. The columns are nearly parallel; this is
why a small Wronskian determinant alone does not produce the cut.

\emph{(b) Two adjacent small rows (orange box).}
At $k=4$, the column norms are approximately
$6.19\times10^{-4}$ and $7.46\times10^{-4}$.
Consequently, $L_4\approx7.46\times10^{-5}<\tau$.
Thus \emph{the orange block gives a valid cut}.
In fact, the adjacent cuts $k=3$ and $k=5$ also satisfy the
condition, owing to their small off-diagonal couplings.

\emph{(c) A large--small--large configuration (green boxes).}
The columns $\gamma_6$ and $\gamma_8$ have norms about $0.7071$,
while $\|\gamma_7\|_2\approx7.46\times10^{-7}$.
Yet $L_6\approx L_7\approx0.1768>\tau$: neither adjacent
cut is valid. This is precisely why the contradiction argument
uses the three-term recurrence to compare the directions
of $\gamma_6$ and $\gamma_8$ directly, skipping the
intermediate small column.

The three cases are therefore not three ways of finding a cut.
They explain how the proof rules out the possibility that
\emph{every} cut has leakage greater than $\tau$.
The corresponding vectors are separated at the same cut:
\begin{equation}\label{eq:localization-vectors-cut}
\resizebox{0.91\linewidth}{!}{$\displaystyle
\begin{pmatrix}v_1^{\mathsf T}\\ v_2^{\mathsf T}\end{pmatrix}
\approx
\begin{tikzpicture}[baseline=(V.center)]
\matrix (V) [matrix of math nodes,
 left delimiter={(},right delimiter={)},
 nodes={inner sep=1.9pt,anchor=center},
 column sep=9pt,row sep=4pt,font=\scriptsize] {
 0.577350 \pgfmatrixnextcell 0.577350 \pgfmatrixnextcell 0.577351 \pgfmatrixnextcell -0.000601 \pgfmatrixnextcell 0 \pgfmatrixnextcell 0 \pgfmatrixnextcell 0 \pgfmatrixnextcell 0 \\
 0 \pgfmatrixnextcell 0 \pgfmatrixnextcell 0 \pgfmatrixnextcell 0 \pgfmatrixnextcell -0.000737 \pgfmatrixnextcell 0.707107 \pgfmatrixnextcell 7.37\cdot10^{-7} \pgfmatrixnextcell -0.707107 \\
};
\coordinate (vCutTop) at ($(V-1-4.north east)!0.5!(V-1-5.north west)$);
\coordinate (vCutBottom) at ($(V-2-4.south east)!0.5!(V-2-5.south west)$);
\coordinate (vCutBottomAligned) at (vCutTop |- V-2-4.south);
\draw[red!80!black,densely dashed,line width=1.25pt]
 ([yshift=7pt]vCutTop) -- ([yshift=-7pt]vCutBottomAligned);
\node[red!80!black,above=5pt,font=\scriptsize] at (vCutTop) {$k=4$};
\end{tikzpicture}.
$}
\end{equation}

\end{example}

\paragraph{Where the structure enters.}
Only the localization lemma uses the tridiagonal recurrence:
through the Wronskian identity, the small-leakage cut, and the
boundary terms created by coordinate truncation. The rank-two
auxiliary construction and the nuclear-norm transfer apply to
general Hermitian matrices whenever appropriate approximate
eigenvectors with separated $X$-responses can be found.

\begin{remark}[Role of the tridiagonal structure]
\label{rem:tridiagonal-localization}
The tridiagonal assumption enters at three points:
the discrete Wronskian identity, the three-term
recurrence used to control directions across an
intermediate small row, and the two boundary terms
created by truncation at a cut. The initial
orthogonality argument and the projection onto
coordinate subspaces are not specific to
tridiagonal matrices.

Thus the structural role of tridiagonality is
to guarantee a cut with small boundary leakage.
Once the approximate eigenvectors in
\Cref{lem:localization} have been constructed,
the auxiliary-matrix argument and the subsequent
slope-variation estimates do not require the
auxiliary matrix to remain tridiagonal.
\end{remark}

\subsection{From the cut to the gap certificate}
\label{subsec:gap-transfer}

The preceding lemma is the only step that uses the tridiagonal
recurrence. We now replace the two localized approximate eigenvectors
by exact eigenvectors with a common eigenvalue in a nearby auxiliary
matrix. This gives a uniform slope response, which we then transfer
back to the original family.

\begin{lemma}[A nearby matrix with large slope variation]
\label{lem:auxiliary}
Under the assumptions of \Cref{lem:localization}, there is a real
symmetric matrix $\widetilde H$ for which
\begin{equation}\label{eq:auxiliary}
\widetilde Hv_i=\frac{\lambda+\mu}{2}v_i\quad(i=1,2),
\qquad\norm{H(t)-\widetilde H}_*\leq96n\sqrt\gamma,
\end{equation}
where $\|\cdot\|_*$ denotes the nuclear norm. 
The nuclear-norm bound measures its distance from
$H(t)$ and will allow us to transfer the slope-variation
estimate back to the original perturbation family. The auxiliary matrix is not required to be tridiagonal.
Define $\widetilde H(s)=\widetilde H+(s-t)X$ and
$F^{\mathrm{aux}}(s)=\sum_{j=1}^n j\lambda_j(\widetilde H(s))$.
Then the centered slope variation of the auxiliary family satisfies
\begin{equation}\label{eq:auxiliary-defect}
\mathcal S_h^{\mathrm{aux}}(t)
:=\frac{F^{\mathrm{aux}}(t+h)+F^{\mathrm{aux}}(t-h)
-2F^{\mathrm{aux}}(t)}{h}\geq1\qquad(h>0).
\end{equation}
\end{lemma}

\begin{proof}
Take $v_1,v_2$ from \Cref{lem:localization} and set
$Q=v_1v_1^{\mathsf T}+v_2v_2^{\mathsf T}$ and
$\widetilde H=(I-Q)H(t)(I-Q)+(\lambda+\mu)Q/2$.
Because $Qv_i=v_i$, both are exact eigenvectors of
$\widetilde H$ at the common eigenvalue $(\lambda+\mu)/2$.
Moreover, expanding $H(t)-\widetilde H$ gives
\[
H(t)-\widetilde H
=\bigl(H(t)-\tfrac{\lambda+\mu}{2}I\bigr)Q
+Q\bigl(H(t)-\tfrac{\lambda+\mu}{2}I\bigr)
-Q\bigl(H(t)-\tfrac{\lambda+\mu}{2}I\bigr)Q.
\]
Since $Q$ has rank two, the nuclear--Frobenius norm inequality
and \eqref{eq:quasi} yield the compact estimate
\[
\norm{H(t)-\widetilde H}_*
\leq3\sqrt2\left(\sum_{i=1}^2
\norm{\bigl(H(t)-\tfrac{\lambda+\mu}{2}I\bigr)v_i}_2^2
\right)^{1/2}\leq96n\sqrt\gamma.
\]

For the slope variation, choose an ordered orthonormal eigenbasis
of $\widetilde H$ in which $v_1,v_2$ occupy adjacent places
$j,j+1$ in their common eigenspace. Assign the corresponding
weights $j,j+1$ first to $v_1,v_2$, then in the opposite order.
Both choices maximize the variational representation
\eqref{eq:variational} at $s=t$. The corresponding matrices
$G_+,G_-$ satisfy
\[
\tr(G_+\widetilde H)=\tr(G_-\widetilde H)=F^{\mathrm{aux}}(t),
\qquad G_+-G_-=v_2v_2^{\mathsf T}-v_1v_1^{\mathsf T}.
\]
Using $G_+$ at $t+h$ and $G_-$ at $t-h$ gives
\[
\mathcal S_h^{\mathrm{aux}}(t)
\geq\tr\bigl((G_+-G_-)X\bigr)
=v_2^{\mathsf T}Xv_2-v_1^{\mathsf T}Xv_1\geq1.
\]
\end{proof}

\begin{proof}[Proof of \Cref{prop:bridge}]
We prove the contrapositive. Suppose
$\gap(H(t))\leq\gamma:=h^2/(2^{20}n^4)$.
Because $0<h\leq1$, we have $0<\gamma\leq1/(256n^2)$.
Choose two eigenvalues $\lambda,\mu$ realizing the minimum gap,
including $\lambda=\mu$ if the gap is zero. By
\Cref{lem:localization,lem:auxiliary}, we have
$\norm{H(t)-\widetilde H}_*\leq96n\sqrt\gamma$ and
$\mathcal S_h^{\mathrm{aux}}(t)\geq1$.

By the variational representation \eqref{eq:variational}
and operator--nuclear norm duality, for every $s$,
\begin{equation}\label{eq:nuclear}
|F(s)-F^{\mathrm{aux}}(s)|
\leq n\norm{H(s)-\widetilde H(s)}_*
=n\norm{H(t)-\widetilde H}_*.
\end{equation}
The original and auxiliary families differ by the constant matrix
$H(t)-\widetilde H$, since
$H(s)-\widetilde H(s)=H(t)-\widetilde H$ for all $s$.
Applying \eqref{eq:nuclear} at $t-h,t,t+h$ yields
\[
\mathcal S_h(t)\geq\mathcal S_h^{\mathrm{aux}}(t)
-\frac{4n}{h}\norm{H(t)-\widetilde H}_*
\geq1-\frac{384n^2\sqrt\gamma}{h}.
\]
At $\gamma=h^2/(2^{20}n^4)$, the right-hand side equals
$1-384/1024=5/8>1/2$, proving \eqref{eq:bridge}.
\end{proof}

\section{The algorithm: adaptive selection of the perturbation parameter}
\label{sec:selection}

The purpose of this section is to construct a diagonal
perturbation $D=tX$ that guarantees quantitative spectral
separation. By \Cref{prop:bridge}, it suffices to find
a parameter $t$ and a window radius $h$ such that
\[
\mathcal S_h(t)<\frac12
\quad\Longrightarrow\quad
\gap(H(t))>\frac{h^2}{2^{20}n^4}.
\]
The key observation is that the convex spectral function
$F$ has a bounded total increase in slope. Although its
slope may vary rapidly over some windows, repeated interval
refinement allows us to find a window with sufficiently
small slope variation.

We construct an interval-selection algorithm that checks the local
certificate at each selected window and terminates as soon as it is
verified. We prove an explicit worst-case iteration bound and then
explain how the required spectral-function evaluations can be carried
out efficiently. A simpler selector, which keeps one quarter of
each window, is analyzed in \Cref{app:basic-selection}; it also
terminates after $O(\log n)$ refinements, but with a weaker polynomial
gap bound. 

\subsection{Adaptive selection by endpoint slope variation}
\label{subsec:improved-selection}

We represent each interval by $I=[t-h,t+h]$, where $t$ is its center
and $h$ is its radius. Since the local gap certificate of
\Cref{prop:bridge} depends on $h^2$, the selector should retain as much
of each window as possible. Comparing slopes near the endpoints allows
us to keep almost half of the parent interval at each iteration.

\paragraph{Endpoint slope variation.}
For a window $I=[t-h,t+h]$ and $0<\theta\leq1/4$, define
\begin{equation}\label{eq:endpoint-slope}
\boxed{
\mathcal S_h^{(\theta)}(t)
:=
\frac{F(t+h)-F(t+(1-2\theta)h)}{2\theta h}
-
\frac{F(t-(1-2\theta)h)-F(t-h)}{2\theta h}.
}
\end{equation}
This measures the difference between the secant slopes near the two
endpoints and is nonnegative by convexity.

The endpoint slopes also control the centered slope variation. Indeed,
for $0<\theta\leq1/2$, convexity gives
\[
\frac{F(t-h+2\theta h)-F(t-h)}{2\theta h}
\leq\frac{F(t)-F(t-h)}h,
\qquad
\frac{F(t+h)-F(t)}h
\leq\frac{F(t+h)-F(t+h-2\theta h)}{2\theta h}.
\]
Subtracting yields
\begin{equation}\label{eq:potential-controls-S}
\mathcal S_h(t)\leq\mathcal S_h^{(\theta)}(t).
\end{equation}
Thus endpoint slope variation can be used for selection, while the
stopping test remains the original condition $\mathcal S_h(t)<1/2$.

Let $K=\lceil\log_2(16n^3)\rceil$ and $\theta=1/(4K)$.
For a current interval $I=[t-h,t+h]$, consider
\[
I_{\mathrm L}=[t-(1-2\theta)h,t],
\qquad
I_{\mathrm R}=[t,t+(1-2\theta)h].
\]
Both candidates have radius $h'=(1/2-\theta)h$.

\begin{algorithm}[t]

\caption{Adaptive selection of a separating perturbation}
\label{alg:selection}
\DontPrintSemicolon
\KwIn{$T$, $n\geq2$, and $0<\delta<1/2$.}
\KwOut{A perturbation parameter $t$ and certification radius $h$.}

Set $K=\lceil\log_2(16n^3)\rceil$ and $\theta=1/(4K)$\;
Set $t=\delta/(2n)$ and $h=\delta/(2n)$\;
Evaluate $F(t-h),F(t),F(t+h)$ to absolute error at most $h/256$\;
\If{$\widehat{\mathcal S}_h(t)+1/64<1/2$}{
    \Return{$t,h$}\;
}

\For{$j=0,\ldots,K-1$}{
    Set $h=(1/2-\theta)h$\;
    Evaluate the values of $F$ required for
    $\mathcal S_h^{(\theta)}(t-h)$ and
    $\mathcal S_h^{(\theta)}(t+h)$, each to absolute
    error at most $2\theta h/256$\;
    Compute $\widehat{\mathcal S}_h^{(\theta)}(t-h)$
    and $\widehat{\mathcal S}_h^{(\theta)}(t+h)$\;
    \eIf{$\widehat{\mathcal S}_h^{(\theta)}(t-h)
       \leq\widehat{\mathcal S}_h^{(\theta)}(t+h)$}{
        Set $t=t-h$\;
    }{
        Set $t=t+h$\;
    }
    Evaluate the additional midpoint value $F(t)$
    to absolute error at most $2\theta h/256$\;
    Compute $\widehat{\mathcal S}_h(t)$\;
    \If{$\widehat{\mathcal S}_h(t)+1/64<1/2$}{
        \Return{$t,h$}\;
    }
}
\end{algorithm}

The following lemma combines the exact contraction of endpoint
slope variation with the effect of approximate evaluations.

\begin{lemma}[Contraction of endpoint slope variation]
\label{lem:potential-contraction}
\label{lem:potential-inexact}
Let $I=[t-h,t+h]$, $h'=(1/2-\theta)h$, and
$0<\theta\leq1/4$. Then
\begin{equation}\label{eq:potential-contraction}
\min\left\{
\mathcal S_{h'}^{(\theta)}(t-h'),
\mathcal S_{h'}^{(\theta)}(t+h')
\right\}
\leq\frac12\mathcal S_h^{(\theta)}(t).
\end{equation}
Moreover, the selected centers and radii of
\Cref{alg:selection} satisfy
\begin{equation}\label{eq:potential-inexact-contraction}
\mathcal S_{h_{j+1}}^{(\theta)}(t_{j+1})
\leq\frac12\mathcal S_{h_j}^{(\theta)}(t_j)+\frac1{32}.
\end{equation}
\end{lemma}

\begin{proof}
Convexity gives a contraction between the parent and its two
candidate windows: the endpoint secants are ordered, and the two child slope differences correspond
to disjoint increments between the parent endpoint slopes. Hence
\[
\mathcal S_{h'}^{(\theta)}(t-h')
+\mathcal S_{h'}^{(\theta)}(t+h')
\leq\mathcal S_h^{(\theta)}(t),
\]
which proves \eqref{eq:potential-contraction}.

For approximate evaluations, each endpoint slope variation
combines four values of $F$ with coefficients $\pm1/(2\theta h_{j+1})$.
The prescribed absolute error $2\theta h_{j+1}/256$
therefore gives an error of at most $4/256=1/64$ in each
computed variation. Selecting the smaller computed value
adds at most $2/64=1/32$, proving
\eqref{eq:potential-inexact-contraction}.
\end{proof}

\begin{theorem}[Constructive spectral separation]
\label{thm:selection}

Algorithm~\ref{alg:selection} terminates after at most $K$ refinements,
where $K=\lceil\log_2(16n^3)\rceil=O(\log n)$,
and returns $t,h$ satisfying
\begin{equation}\label{eq:selection-output}
\mathcal S_h(t)<\frac12,
\qquad h>\frac{\delta}{128n^4}.
\end{equation}
Consequently, $D=tX$ satisfies
\begin{equation}\label{eq:selection-gap}
\boxed{
\norm D_2\leq\delta,\qquad
\gap(T+D)>2^{-34}\frac{\delta^2}{n^{12}}.
}
\end{equation}
The algorithm requires at most $3+8K$ spectral-function evaluations.
\end{theorem}

\begin{proof}
Initially $t_0=h_0=\delta/(2n)$, so the initial interval is
$[0,\delta/n]$, and each refinement selects a subinterval of its
parent. Hence every selected window lies in $[0,\delta/n]$ and
$\norm{tX}_2=nt\leq\delta$.

By \Cref{lem:convex} and
\eqref{eq:potential-inexact-contraction},
\[
\mathcal S_{h_K}^{(\theta)}(t_K)
\leq 2^{-K}n^3+\frac1{32}\sum_{j=0}^{K-1}2^{-j}
<\frac18.
\]
By \eqref{eq:potential-controls-S}, we also have
$\mathcal S_{h_K}(t_K)<1/8$. The endpoint values and
the additional midpoint evaluation give
\[
\left|\widehat{\mathcal S}_{h_K}(t_K)
-\mathcal S_{h_K}(t_K)\right|
\leq\frac{\theta}{32}\leq\frac1{128}.
\]
Thus $\widehat{\mathcal S}_{h_K}(t_K)+1/64<1/2$,
so the certified stopping test succeeds by iteration $K$.

It remains to bound the radius. Termination at an earlier
iteration only increases the returned radius. Since
\[
L_K=\frac{\delta}{n}
\left(\frac12-\theta\right)^K
=\frac{\delta}{n}2^{-K}
\left(1-\frac1{2K}\right)^K
>\frac{\delta}{64n^4},
\]
where we used Bernoulli's inequality and $2^K<32n^3$,
the returned radius satisfies
$h\geq L_K/2>\delta/(128n^4)$.
Applying \Cref{prop:bridge} yields
\[
\gap(T+tX)>
\frac{h^2}{2^{20}n^4}
>2^{-34}\frac{\delta^2}{n^{12}}.
\]
Finally, the initial test uses three evaluations of $F$,
and each refinement uses at most eight.
\end{proof}

\subsection{Cost of the spectral-function evaluations}
\label{subsec:oracle}

For the selector, the required absolute accuracy in each endpoint
comparison is $2\theta h/256$, where $h$ is the updated candidate
radius; the additional centered-slope test uses the same accuracy.
Since the final search radius after $K$ refinements is greater than
$\delta/(128n^4)$, the smallest requested accuracy satisfies
\begin{equation}\label{eq:minimum-oracle-accuracy}
\varepsilon_{\min}>
\frac{\delta}{2^{16}Kn^4},
\qquad\text{so}\qquad
\log(n/\varepsilon_{\min})=O(1+\log(n/\delta)).
\end{equation}
Thus \Cref{alg:selection} requires at most $3+8K=O(\log n)$ values of
$F$, each to an absolute accuracy no finer than $\varepsilon_{\min}$.
If approximations $\widehat\lambda_j(t)$ satisfy
$|\widehat\lambda_j(t)-\lambda_j(t)|\leq\eta$, then
\[
\Bigl|\sum_{j=1}^n j\widehat\lambda_j(t)-F(t)\Bigr|
\leq\frac{n(n+1)}2\eta,
\]
so $\eta=2\varepsilon/[n(n+1)]$ yields accuracy $\varepsilon$ for $F$.
The cost of the construction is therefore that of $O(\log n)$
computations of all eigenvalues of a Hermitian tridiagonal matrix to an
absolute accuracy that is polynomial in $\delta/n$; no eigenvectors are
required.

Any eigenvalue solver with this accuracy guarantee can be used.
Sobczyk~\cite{Sobczyk2025} gives a complete analysis, in the
Real RAM model, of the divide-and-conquer eigensolver of Gu and
Eisenstat~\cite{GuEisenstat1995} accelerated by the fast multipole
method, with the required guarantees of the fast multipole routine
established in that work \cite[Proposition~1.1]{Sobczyk2025}. In
particular, all eigenvalues of an unreduced real symmetric tridiagonal
matrix with spectral norm at most one can be approximated to absolute
accuracy $\eta$ using $O(n\log^{c}(n/\eta))$ arithmetic operations, for
an absolute constant $c$ \cite[Appendix~C.2, Corollary~C.1 of the full
version]{Sobczyk2025}. This applies to $H(t)$ after the diagonal unitary
similarity that makes its off-diagonal entries real and nonnegative
(which commutes with $tX$), after scaling by $(1+\delta)^{-1}$, and,
when some off-diagonal entry vanishes, to each unreduced diagonal block
followed by a merge of the sorted eigenvalue lists. We do not
pursue the complexity of this step further.

The stopping and selection tests above assume
the stated absolute evaluation-error guarantees.
An implementation using ordinary floating-point
approximations without such guarantees may still
produce a valid perturbation whenever its final
centered-slope test is rigorously certified, but
the worst-case iteration bound requires the
endpoint-selection errors to be controlled.
These statements concern exact evaluations with
prescribed absolute errors; they do not address
finite-precision or bit complexity.

\paragraph{Related work: eigenvalue algorithms and complexity.}
Since the selector reduces the construction to $O(\log n)$
spectral-function evaluations, its overall cost rests on how
efficiently all eigenvalues of a tridiagonal matrix can be
approximated.

Classical structured eigensolvers, including
the divide-and-conquer algorithm of Gu and
Eisenstat~\cite{GuEisenstat1995}, provide
important computational foundations.
Sobczyk~\cite{Sobczyk2025} gives deterministic
complexity bounds for Hermitian eigenproblems,
including a nearly quadratic full diagonalization and a nearly
linear eigenvalues-only computation for symmetric tridiagonal matrices,
both with complete proofs in the Real RAM model, together with
deterministic subroutines for spectral gaps and spectral projectors.
Such results motivate the use of a spectral
oracle that computes the weighted eigenvalue
sum without requiring a complete eigenvector
matrix. Question~4 of \cite{Sobczyk2025}, also discussed at the
2025 Householder Symposium, is motivated by the observation that the
perturbations used by Gu and Eisenstat for arrowhead matrices produce
precisely the required eigenvalue separation, whereas transferring such a
perturbation from the arrowhead level back to the original tridiagonal (or
Hermitian) matrix involves nontrivial steps. Our perturbation is instead
constructed directly in the original coordinates of $T$, so no such
transfer is required.

Related algorithmic developments include the
randomized subspace methods of Nakatsukasa and
Tropp~\cite{NakatsukasaTropp2024}, the contour
and subspace approaches studied by Horning and
Townsend~\cite{HorningTownsend2020}, and rigorous
spectral computations developed by Colbrook
and Townsend~\cite{ColbrookTownsend2024,
ColbrookTownsend2025}. These works address
different eigenvalue-computation problems, but
highlight the importance of scalable algorithms
and reliable spectral information.

Our contribution is complementary: we construct
a small deterministic perturbation that gives
a guaranteed spectral gap, while explicitly
accounting for the cost of selecting and
certifying that perturbation.

\section{Numerical experiments}
\label{sec:numerics}

The theoretical guarantee of \Cref{thm:selection} is 
conservative. We investigate how the selector of \Cref{alg:selection} behaves in
practice, whether its centered-slope stopping test can be certified
in finite precision, and how the resulting spectral gaps compare
with their analytical lower bounds.

Throughout, we use $\delta=0.48$ and the fixed perturbation direction
$X=\diag(1,\ldots,n)$. The experiments distinguish three quantities:
the worst-case iteration bound $K$, the actual termination iteration
$j_*$, and the certified minimum gap of the returned matrix.

\paragraph{Test matrices and implementation.}
\label{subsec:numerical-setup}
We use the six families of real symmetric tridiagonal matrices in
\Cref{tab:test-families} (\Cref{app:test-families}), normalized so
that $\|T\|_2\leq1$; by \Cref{rem:real}, separate complex Hermitian tests
are not needed. Double-precision tridiagonal eigensolvers propose the
successive windows, and each proposed window is tested with rigorous
enclosures of the centered slope variation. These come from a structured
oracle: double-precision eigenvalues give initial locations, individual
eigenvalues are enclosed by two-sided Sturm bracketing validated against
the input matrix in interval arithmetic, unresolved clusters are enclosed
and counted collectively, and the precision is raised only when an
enclosure cannot decide the stopping test. The ordered enclosures give
rigorous bounds on $F(t)$, $\mathcal S_h(t)$ and the minimum gap. This
oracle has quadratic cost in $n$; it demonstrates rigorous
certification, not the cost of the eigenvalue solvers discussed in
\Cref{subsec:oracle}. We report the first
window along the double-precision selection path whose centered slope
variation is certified below $1/2$; a separate rigorous fallback controls
the endpoint comparisons required by the worst-case termination proof.

\begin{figure}[t]
\centering
\includegraphics[width=\linewidth]{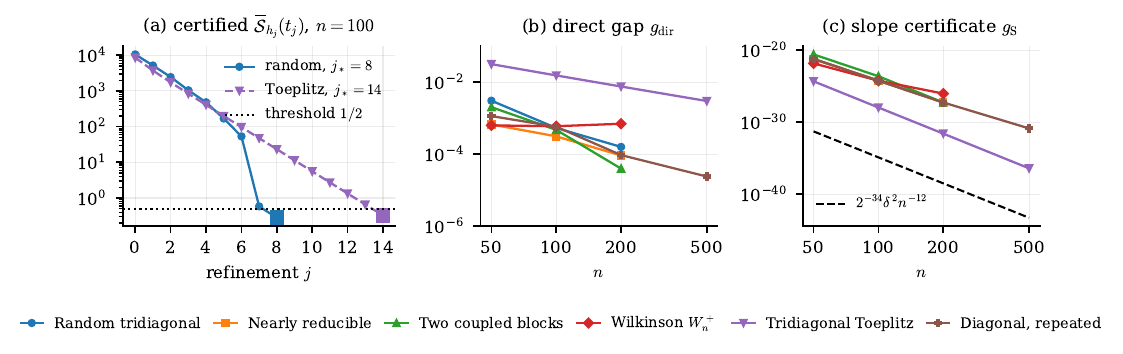}
\caption{Certified numerical behavior (strict-Sturm dataset, seed~1).
(a) Rigorous upper bound on the centered slope variation along the
selection path at $n=100$; squares mark the first certified stop $j_*$.
(b) Directly certified minimum gap $g_{\mathrm{dir}}$ of the returned
matrix. (c) Slope-based certificate $g_{\mathrm S}$ at the stopping window,
with the universal lower bound $2^{-34}\delta^2n^{-12}$ (dashed). The
order $n=500$ was run for two families.}
\label{fig:numerics}
\label{fig:certified-iterations}
\label{fig:certified-gaps}
\end{figure}

\paragraph{Iteration counts.}
\label{subsec:numerical-iterations}
\Cref{alg:selection} has the worst-case bound
$K=\lceil\log_2(16n^3)\rceil$, but the certified stopping test
$\overline{\mathcal S}_{h_{j_*}}(t_{j_*})<1/2$ is usually met much
earlier (\Cref{tab:certified-results}). For a random tridiagonal matrix
of order $100$ the algorithm stops after $j_*=8$ refinements, compared with
$K=24$, and for the Toeplitz matrix of the same order after $14$; the
Toeplitz family generally needs the most refinements. \Cref{fig:numerics}(a)
shows the mechanism. Along the Toeplitz path the certified slope variation
decreases by a factor of about two per refinement, the rate guaranteed for
the endpoint variation by \Cref{lem:potential-contraction}, whereas along
the random path it drops by two orders of magnitude in a single step once
the window avoids the nearby collisions. These observations concern the tested instances and do not
improve the worst-case bound.

\paragraph{Analytical versus directly certified gaps.}
\label{subsec:numerical-certificates}
For the returned parameter and radius $(t_*,h_*)$,
\Cref{prop:bridge} and the certified eigenvalue enclosures
$\lambda_j(T+t_*X)\in[\underline\lambda_j,\overline\lambda_j]$ give the
two lower bounds
\[
g_{\mathrm S}:=\frac{h_*^2}{2^{20}n^4}<\gap(T+t_*X),
\qquad
g_{\mathrm{dir}}:=\min_{1\leq j<n}
(\underline\lambda_{j+1}-\overline\lambda_j)\leq\gap(T+t_*X),
\]
the second of which is not required by the selection theorem.
\Cref{fig:numerics}(b)--(c) compares them with the universal guarantee
$2^{-34}\delta^2n^{-12}$. The perturbation is effective where it is needed most: the
repeated diagonal matrix has a multiple eigenvalue and $W_{200}^+$ has
$\gap(T)\approx4\times10^{-317}$, while every returned matrix in
\Cref{tab:certified-results} has a certified gap above $10^{-5}$. For the
random matrix with $n=100$, $g_{\mathrm S}=5.98\times10^{-25}$ while
$g_{\mathrm{dir}}=5.31\times10^{-4}$, more than twenty orders of magnitude
larger. Two sources of conservatism are visible: early termination
($j_*<K$) retains a larger window and hence a stronger $g_{\mathrm S}$,
but even then the universal local certificate lies far below the actual
gap, most markedly for the Toeplitz family. This identifies the
localization-to-gap estimate of \Cref{sec:local-gap} as the natural
target for improvement, without by itself implying a better worst-case
exponent than $n^{-12}$.

\begin{table}[t]
\centering
\caption{Representative runs of the strict-Sturm implementation
with $\delta=0.48$. The index $j_*$ is the first rigorously
certified stopping iteration; $\gap(T)$ is the minimum gap
before perturbation (for $W_{200}^+$ computed by Sturm bisection in
420-digit arithmetic, otherwise in double precision).}
\label{tab:certified-results}
\small
\begin{tabular}{lrrrrr}
\hline
Family and $n$ & $K$ & $j_*$ & $\gap(T)$ & $g_{\mathrm{dir}}$ & $g_{\mathrm S}$\\
\hline
Random, 100       &24& 8&$1.67\times10^{-4}$&$5.31\times10^{-4}$&$5.98\times10^{-25}$\\
Toeplitz, 100     &24&14&$1.45\times10^{-3}$&$1.54\times10^{-2}$&$1.13\times10^{-28}$\\
Random, 200       &27&10&$9.09\times10^{-6}$&$1.58\times10^{-4}$&$5.63\times10^{-28}$\\
Wilkinson, 200    &27& 8&$4.34\times10^{-317}$&$6.97\times10^{-4}$&$9.71\times10^{-27}$\\
Toeplitz, 200     &27&17&$3.66\times10^{-4}$&$7.57\times10^{-3}$&$2.65\times10^{-32}$\\
Diagonal, 500     &31&12&$0$&$2.37\times10^{-5}$&$1.42\times10^{-31}$\\
Toeplitz, 500     &31&21&$5.90\times10^{-5}$&$2.99\times10^{-3}$&$4.04\times10^{-37}$\\
\hline
\end{tabular}
\end{table}

\section{Conclusions and outlook}
\label{sec:conclusion}

We have given a deterministic construction of a small diagonal
perturbation that produces inverse-polynomial spectral separation for
arbitrary Hermitian tridiagonal matrices. The proof combines
tridiagonal invariant-subspace localization with the variational geometry
of a convex weighted spectral function: its finite-scale slope variation
provides a local gap certificate, and an adaptive, error-tolerant window
selector finds a certifying perturbation after $O(\log n)$ refinements.
The construction needs only $O(\log n)$ evaluations of a weighted
eigenvalue sum to inverse-polynomial accuracy, so its cost is governed
by that of tridiagonal eigenvalue computations.
Rigorous experiments show that the stopping test often succeeds well
before the theoretical bound, and that the directly certified gaps are
far larger than the universal certificate; the localization-to-gap
constants and the dimension dependence are therefore natural targets for
a sharper analysis.

Further questions include the arithmetic and bit complexity of
the spectral evaluations, extensions beyond the tridiagonal class, and the use of
diagonal separation as a preprocessing step in eigensolvers
\cite{NakatsukasaTropp2024,ColbrookTownsend2024,GrafTownsend2026}.
The guiding principle of the construction, that a verifiable structural
certificate replaces an expensive global computation, also underlies our
earlier work on structured recurrence matrices and polynomial
certificates \cite{ZhuNakatsukasa2026VA,ZhuNie2026ToeplitzBW,CartisZhu2026SoSExact}.

\paragraph{Deterministic pseudospectral shattering.}
Sobczyk's Question~4 \cite{Sobczyk2025} asks for deterministic
pseudospectral shattering for general Hermitian and even non-Hermitian
matrices. The construction here settles the eigenvalue-gap version for
Hermitian tridiagonal matrices with a diagonal perturbation. Whether it
can be extended to dense Hermitian or non-Hermitian matrices, possibly
in combination with the arrowhead perturbations of Gu and Eisenstat
\cite{GuEisenstat1995} and the deterministic diagonalization framework
of \cite{Sobczyk2025}, is an interesting open direction.

\needspace{8\baselineskip}
\section*{Statements and Declarations}
\paragraph{Acknowledgements.}
We thank Aleksandros Sobczyk for helpful comments on the history of the
deterministic perturbation problem and on the fast multipole results
of~\cite{Sobczyk2025}.

\paragraph{Funding.}
The research of Wenqi Zhu was supported by the Engineering and Physical
Sciences Research Council under grant EP/Y028872/1, \emph{Mathematical
Foundations of Intelligence: An ``Erlangen Programme'' for AI}.
Ping Nie gratefully acknowledges the support of AI4ALPHA and NetMind.AI.

\paragraph{Competing interests.}
The authors have no competing interests to declare that are relevant to
the content of this article.

\paragraph{Data availability.}
The certified implementation, the six families of test matrices and the
scripts producing the figures in Section~\ref{sec:numerics} are
publicly available at\newline
\mbox{\url{https://github.com/KateWenqiZhu/spectral-separation-certificates}}.

\paragraph{Use of generative AI and AI-assisted technologies.}
During the preparation of this work, the authors used AI tools (Claude, Anthropic; Codex, OpenAI; Qwen App, Alibaba) to generate candidate constructions and code, perform computations, draft mathematical arguments, and improve the readability of the text. After using these
tools, the authors reviewed and edited the content as needed and take full
responsibility for the content of the publication. The theorems are
established by the written proofs; the numerical experiments use rigorous
Sturm-based eigenvalue enclosures to certify the reported quantities.

\appendix

\section{A basic selector}
\label{app:basic-selection}
\label{subsec:basic-selection}

We represent each interval by $I=[t-h,t+h]$,
where $t$ is its center and $h$ is its radius.
At each iteration, we compare the centered slope
variations on the two middle subintervals
\[
I_{\mathrm L}=[t-h/2,t],
\qquad
I_{\mathrm R}=[t,t+h/2].
\]
Their centers are $t-h/4$ and $t+h/4$, respectively,
and both have radius $h/4$.

The basic procedure selects the subinterval with
smaller computed slope variation. At each iteration,
we also check whether the selected window already
satisfies the local gap certificate. We will show
that the algorithm must terminate after at most
$O(\log n)$ refinements, while maintaining an
inverse-polynomial lower bound on the window radius.

\begin{algorithm}[H]

\caption{Basic selection of a perturbation parameter}
\label{alg:selection-basic}
\DontPrintSemicolon
\KwIn{$T$, $n\geq2$, and $0<\delta<1/2$.}
\KwOut{A perturbation parameter $t$ and radius $h$.}

Set $t=\delta/(2n)$ and $h=\delta/(2n)$\;
Let $K_{\mathrm b}$ be the smallest nonnegative integer
such that $(3/2)^{K_{\mathrm b}}\geq16n^3$\;

Evaluate $F(t-h)$, $F(t)$, and $F(t+h)$
to absolute error at most $h/256$\;
\If{$\widehat{\mathcal S}_h(t)+1/64<1/2$}{
    \Return{$t,h$}\;
}

\For{$j=0,\ldots,K_{\mathrm b}-1$}{
    Set $h=h/4$\;
    Evaluate $F$ at $t-2h,t-h,t,t+h,t+2h$
    to absolute error at most $h/256$\;
    Compute $\widehat{\mathcal S}_h(t-h)$ and
    $\widehat{\mathcal S}_h(t+h)$\;
    \eIf{$\widehat{\mathcal S}_h(t-h)
    \leq\widehat{\mathcal S}_h(t+h)$}{
        Set $t=t-h$\;
    }{
        Set $t=t+h$\;
    }
    \If{$\widehat{\mathcal S}_h(t)+1/64<1/2$}{
        \Return{$t,h$}\;
    }
}
\end{algorithm}

\begin{example}[One step of interval selection]
\label{ex:selection}
Consider the spectral function in \Cref{ex:collision},
whose centered slope variation is
\[
\mathcal S_h(t)
=\max\left\{1-\frac{|t-0.1|}{h},0\right\}.
\]
Starting with the illustrative window
$[0.06,0.14]$, centered at $t=0.1$ with $h=0.04$,
we have $\mathcal S_{0.04}(0.1)=1$.

After one refinement, the two candidate windows
are $[0.08,0.10]$ and $[0.10,0.12]$, with centers
$0.09$ and $0.11$ and radius $0.01$. Since $F$
is affine on each window,
\[
\mathcal S_{0.01}(0.09)
=\mathcal S_{0.01}(0.11)=0.
\]
Thus either candidate already satisfies
$\mathcal S_h(t)<1/2$.
More generally, the algorithm selects the candidate
with smaller slope variation, moving towards a
window where the spectral response is more stable.
\end{example}

We next combine the exact contraction property
of centered slope variation with the effect of
approximate spectral-function evaluations.

\begin{lemma}[Contraction of centered slope variation]
\label{lem:defect-contraction}
\label{lem:inexact-contraction}
For a parent interval $I=[t-h,t+h]$, let
$I_{\mathrm L}=[t-h/2,t]$ and
$I_{\mathrm R}=[t,t+h/2]$ be its two middle
subintervals. Then
\begin{equation}\label{eq:defect-contraction}
\min\left\{
\mathcal S_{h/4}\left(t-\frac h4\right),
\mathcal S_{h/4}\left(t+\frac h4\right)
\right\}
\leq\frac23\mathcal S_h(t).
\end{equation}
Moreover, if each value of $F$ is evaluated with
absolute error at most one 256th of the child
radius, the selected center and radius after
$j$ iterations of \Cref{alg:selection-basic}
satisfy
\begin{equation}\label{eq:inexact-contraction}
\mathcal S_{h_{j+1}}(t_{j+1})
\leq
\frac23\mathcal S_{h_j}(t_j)+\frac1{32}.
\end{equation}
\end{lemma}

\begin{proof}
Recall from \eqref{eq:slope-curvature} that
the centered slope variation is a triangularly
weighted integral of $F''$ whenever $F$ is
twice continuously differentiable. For general
convex $F$, the same representation holds with
its distributional second derivative interpreted
as a nonnegative measure $\nu$.

Thus
\[
\mathcal S_h(t)
=\int_{\mathbb R}w_{t,h}(x)\,d\nu(x),
\qquad
w_{t,h}(x)=
\left(1-\frac{|x-t|}{h}\right)_+.
\]
A direct comparison of the triangular weights gives
\[
w_{t-h/4,h/4}(x)+w_{t+h/4,h/4}(x)
\leq\frac43w_{t,h}(x).
\]
Integrating against the nonnegative measure $\nu$
and taking half the sum proves
\eqref{eq:defect-contraction}.

For approximate evaluations, let $r=h/4$
be the child radius. Since the centered slope
variation uses the coefficients $1,-2,1$,
evaluating each of its three function values
within $r/256$ gives
\[
\left|
\widehat{\mathcal S}_r(s)
-\mathcal S_r(s)
\right|
\leq
\frac{(1+2+1)(r/256)}r
=\frac1{64}.
\]
Selecting the smaller computed child variation
therefore incurs an additional error of at most
$2/64=1/32$. Combining this estimate with
\eqref{eq:defect-contraction} proves
\eqref{eq:inexact-contraction}.
\end{proof}

\begin{figure}[htbp]
\centering
\includegraphics[width=0.5\linewidth]{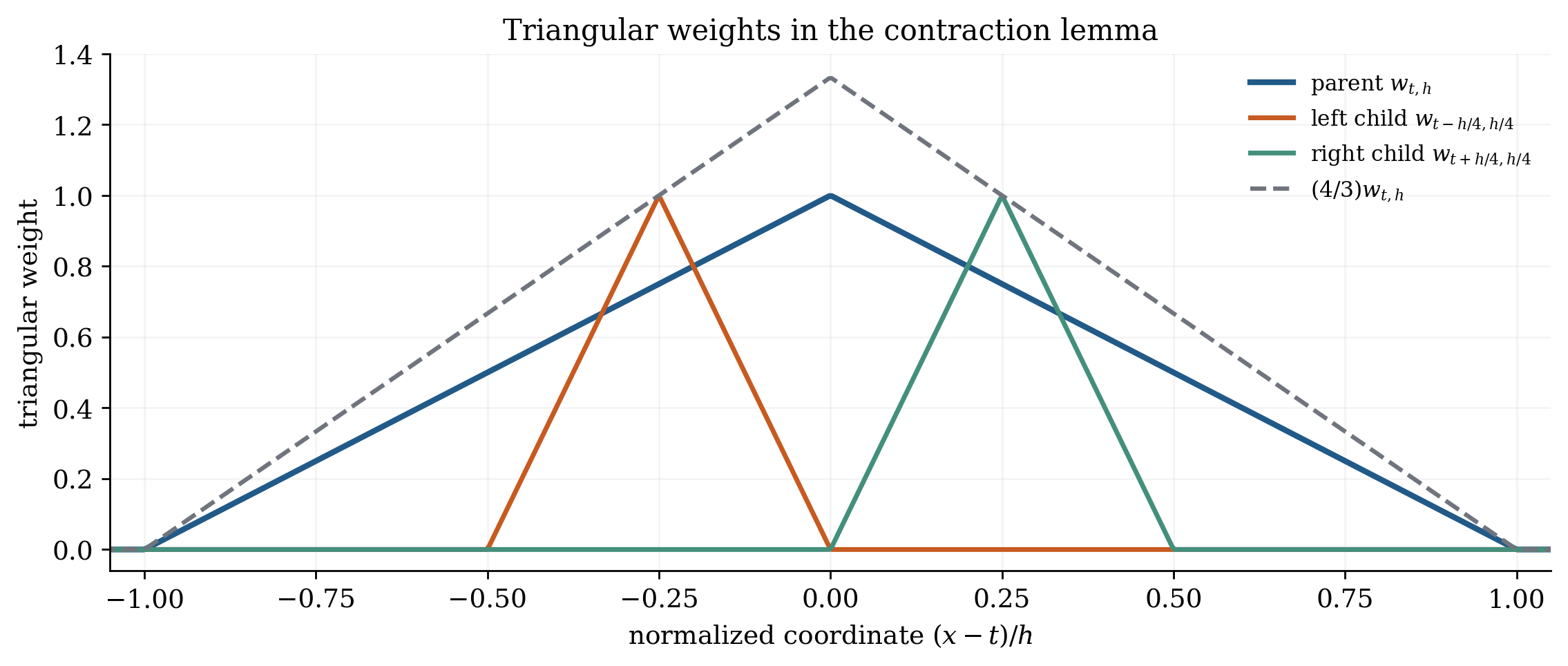}
\caption{Triangular weights for
\Cref{lem:defect-contraction}, plotted in the normalized
coordinate $(x-t)/h$. The two child weights correspond
to the middle subintervals $[t-h/2,t]$ and
$[t,t+h/2]$. Their sum is bounded pointwise by
$(4/3)$ times the parent weight.}
\label{fig:triangular-weights}
\end{figure}

The contraction estimate guarantees that the
slope variation eventually becomes small enough
for the certified stopping test. We now bound
the number of required refinements.

\begin{theorem}[Basic iteration bound and certified scale]
\label{thm:selection-basic}
Algorithm~\ref{alg:selection-basic} terminates after
at most
\[
K_{\mathrm b}
=\left\lceil\log_{3/2}(16n^3)\right\rceil
=O(\log n)
\]
refinements and returns $t,h$ satisfying
\begin{equation}\label{eq:selection-basic-output}
[t-h,t+h]\subseteq[0,\delta/n],
\qquad
\mathcal S_h(t)<\frac12,
\qquad
h\geq\frac{\delta}{2n}4^{-K_{\mathrm b}}.
\end{equation}
It requires at most $3+5K_{\mathrm b}$
spectral-function evaluations.

Consequently, the returned perturbation $D=tX$
satisfies $\norm D_2\leq\delta$ and
\[
\gap(T+D)>\frac{h^2}{2^{20}n^4}.
\]
\end{theorem}

\begin{proof}
Initially $t_0=h_0=\delta/(2n)$, so the initial
interval is $[0,\delta/n]$. Each refinement
selects a subinterval contained in its parent,
with radius one quarter of the parent radius.
Thus, if the algorithm stops after $j$
refinements,
\[
h_j=\frac{\delta}{2n}4^{-j},
\qquad
0\leq j\leq K_{\mathrm b}.
\]

By the secant-slope bounds in \Cref{lem:convex},
the initial slope variation satisfies
$\mathcal S_{h_0}(t_0)<n^3$.
Iterating \eqref{eq:inexact-contraction} gives
\[
\mathcal S_{h_j}(t_j)
<
\left(\frac23\right)^j n^3
+\frac1{32}\sum_{\ell=0}^{j-1}
\left(\frac23\right)^\ell
<
\left(\frac23\right)^j n^3+\frac3{32}.
\]
At $j=K_{\mathrm b}$, the first term is at
most $1/16$, so
\[
\mathcal S_{h_{K_{\mathrm b}}}
(t_{K_{\mathrm b}})
<\frac5{32}.
\]

Every computed slope variation has absolute
error at most $1/64$. Consequently, at the
$K_{\mathrm b}$th refinement,
\[
\widehat{\mathcal S}_{h_{K_{\mathrm b}}}
(t_{K_{\mathrm b}})
+\frac1{64}
<
\frac5{32}+\frac2{64}
<\frac12.
\]
Hence the certified stopping test must succeed
by this iteration. If it succeeds earlier,
the returned window is larger.

The initial test uses three evaluations of $F$,
and each refinement uses at most five.
Finally, $\norm{tX}_2=nt\leq\delta$ and
\Cref{prop:bridge} gives the claimed gap bound.
\end{proof}

\begin{remark}[A lower bound on the window radius]
\label{rem:basic-selection-scale}
Although the algorithm repeatedly reduces the window
radius, it performs only $O(\log n)$ refinements.
Let $\alpha=\log_{3/2}4$. By the minimality of
$K_{\mathrm b}$,
\begin{equation}\label{eq:basic-selection-scale}
h>\frac{\delta}{8n}(16n^3)^{-\alpha}.
\end{equation}
Thus the window radius remains bounded below by
an inverse polynomial in $n$. \Cref{alg:selection}
improves this bound by retaining longer intervals
at each refinement.
\end{remark}

\begin{remark}[Improvement over the basic selector]
\label{rem:selection-scale}
Compared with \Cref{alg:selection-basic}, the selector of
\Cref{alg:selection} achieves a contraction factor of $1/2$ while retaining almost
half of the parent interval. Both algorithms use the same
centered-slope stopping test and the same local gap certificate.
The larger retained radius gives the stronger $n^{-12}$ bound
in \eqref{eq:selection-gap}, without changing the localization
argument of \Cref{sec:local-gap}.
\end{remark}

The basic selector requires at most
$O(\log n)$ refinements and inverse-polynomial
evaluation accuracy, and hence also $O(\log n)$
spectral-function evaluations to inverse-polynomial
accuracy.

\section{Test matrix families}
\label{app:test-families}

\Cref{tab:test-families} lists the six families of real symmetric
tridiagonal matrices used in \Cref{sec:numerics}. 

\begin{table}[H]
\centering\small
\caption{Test families and the spectral features they probe.}
\label{tab:test-families}
\begin{tabular}{p{0.22\linewidth}p{0.43\linewidth}p{0.25\linewidth}}
\hline
Family & Construction & Main purpose\\
\hline
Random tridiagonal & Independent random diagonal and first off-diagonal entries.
& Typical inputs.\\
Nearly reducible & Random tridiagonal matrices with selected first
 off-diagonal entries reduced to $10^{-9}$. & Weakly coupled blocks.\\
Two coupled blocks & Two identical tridiagonal blocks connected by a
 $10^{-7}$ edge. & Nearly repeated pairs.\\
Wilkinson $W_n^+$ & Diagonal entries $|i-(n+1)/2|$ and unit off-diagonals,
 followed by scaling. & Extremely tight clusters.\\
Toeplitz & Zero diagonal and constant off-diagonal entries $1/2$.
 & Known structured spectrum.\\
Repeated diagonal & Diagonal matrices with deliberately repeated entries.
 & Exact multiplicity.\\
\hline
\end{tabular}
\end{table}

\scriptsize{
\bibliographystyle{plain}
\bibliography{sample}
}
\end{document}